\documentclass[a4paper]{article}
\usepackage{amsmath}
\usepackage{amssymb}
\usepackage[english]{babel}
\usepackage{csquotes}
\usepackage{graphicx}
\usepackage{pgfplots}
\usepackage{float}
\usepackage{caption}
\usepackage{hyperref}
\usepackage[left=2.5cm, right=2.5cm, top=2.5cm, bottom=2.5cm]{geometry}
\usepackage{pifont}
\usepackage{tabularx}
\pgfplotsset{compat=1.18}
\usepackage{xurl}
\usepackage[maxbibnames = 10, giveninits = true]{biblatex}
\usetikzlibrary{arrows.meta,calc,fit}

\title{Global Optimization of Beam Profiles for Suppressing the Bump-On-Tail Instability}

\author{Sandro Mignelli, Lukas Einkemmer}
\date{\today}

\begin{document}

\maketitle

\begin{abstract}
    A particle beam with constant intensity that propagates into a plasma triggers a bump-on-tail instability. It has been demonstrated that beam modulation can be used to significantly reduce this instability. However, finding the optimal beam profiles that accomplish this is a complex optimization problem that suffers from an extremely rugged optimization landscape and a relatively high-dimensional parameter space. In this paper, we evaluate a number of optimization algorithms and find that population based methods (in particular, CMA-ES and Differential Evolution) perform best. In particular, we investigate beam profiles where the beam is either completely turned ON or OFF.~This configuration is easier to accomplish experimentally, but leads to a higher-dimensional and thus more challenging optimization problem. We obtain ON/OFF beam profiles that reduce the instability by more than 80\% and thus perform significantly better than what has been reported in the literature. 
\end{abstract}

\section{Introduction}
    A fast particle beam that propagates into a stationary plasma creates a peak in the velocity distribution. This triggers the so-called bump-on-tail instability which converts kinetic energy into electric energy and creates significant density variations in the plasma. In technological applications this is generally undesirable. For example, the bump-on-tail instability created by neutral beam heating in magnetic confinement fusion devices (such as tokamaks) can degrade heating efficiency.

    The bump-on-tail instability has been studied extensively, see, e.g.~\cite{Denavit1985,Carlevaro2015,Zeeland2021}, and it has been demonstrated that modulated beams can have a stabilizing effect. We mention~\cite{Nakamura1970} where it was experimentally demonstrated that a sinusoidally modulated charged particle beam significantly reduces a two-stream instability. More recently, in~\cite{Zeeland2021} neutral beam heating was performed in the DIII-D tokamak, and it has been shown that switching the beam ON/OFF at different frequencies significantly alters the Alfvén mode spectrum. From a theoretical point of view, in~\cite{Qin2014} a linear analysis has shown that modulating the velocity of the beam can have a stabilizing effect. Later in~\cite{10.1063/5.0238408} computer simulations and optimization were used to find periodic beam profiles that reduce the bump-on-tail instability.

    Recently, there has also been some interest from the mathematics community to use PDE-constraint optimization to study the control of a range of kinetic plasma systems. These have mostly focused on finding external electric or magnetic fields to confine, shape, or stabilize plasma systems (as opposed to beam-plasma configurations). Approaches that obtain the control law from theory~\cite{Glass_2012,lukascontrolofinstability,Albi_2025,Chen_2025,albi2025robust}, gradient based optimization~\cite{Bartsch_2024,Crouseilles_2025,Guerra_2026,lu2025controlling}, and a combination of both global and local optimization~\cite{lukassuppressinginstability} have been proposed. These works can also be roughly divided into schemes that require continuous feedback~\cite{Bartsch_2024,Albi_2025,albi2025robust} and systems that require some form of limited feedback~\cite{lukascontrolofinstability,Crouseilles_2025,Chen_2025,lu2025controlling} (most commonly the knowledge of the perturbation at the onset of the instability or information on the moments). 

    A particularly interesting aspect of beam-plasma instabilities, in our view, is that they do not necessarily require feedback, which for micro-instabilities is difficult to obtain, both because the dynamics is very fast and because having detailed measurements of the velocity distribution is difficult. In this paper, we study the bump-on-tail instability introduced by a charged particle beam propagating into a one-dimensional electrostatic plasma. This is modelled by the Vlasov--Poisson equation as has been done in~\cite{10.1063/5.0238408}. In that work, it was observed that the optimization landscape for realistic configurations (i.e.~device size large compared to the Debye length) is very rugged. Our first goal in this paper is thus to systematically study what optimization approach gives good results for this problem. We find that population based approaches perform best and that gradient information can be extremely misleading and, in many cases, deteriorates the performance of the optimization algorithm. In addition, methods based on surrogate models (such as Bayesian optimization) can struggle with the relatively high-dimensional parameter space.

    In~\cite{10.1063/5.0238408} a periodic beam profile has been investigated, and it was shown that this can drastically reduce the severity of the instability. The assumption of periodicity is convenient as it reduces the number of parameters required, i.e.~one does not have to optimize over the entire duration the beam is switched ON.~In this paper, we also investigate the configuration of a beam that is either completely switched ON or OFF (we call this a binary beam). This is a beam shape that has been demonstrated experimentally for beam heating in~\cite{Zeeland2021}. We thus consider it easier to achieve experimentally. However, it also results in a higher-dimensional and thus more challenging optimization problem. In this paper, we demonstrate that binary beam profiles can be found that reduce the instability by more than 80\%, improving upon the result reported in~\cite{10.1063/5.0238408} for a periodically repeating beam-intensity profile.

    The structure of this paper is as follows: In section 2 we introduce the kinetic model and describe the optimization problem. In section~\ref{sec:comparison} we introduce the different optimization algorithms that we investigated and analyze their performance. In section 4 we investigate the binary beam profiles and the associated optimization problem. Finally, in section 5 we summarize our findings and give an outlook on future work.

\section{Problem Statement\label{sec:problem}}
\begin{figure}[H]
        \centering
        \includegraphics[width = 1.0\linewidth]{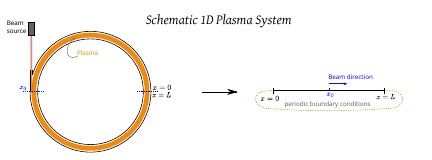}
        \caption{This schematic illustrates the 1+1D Vlasov--Poisson plasma system used in this study. The beam is injected into the plasma at $x_0$. Periodic boundary conditions are applied in the spatial direction, allowing particles and waves to re-enter the domain and enabling the study of plasma interacting with the beam. This illustration is a modification of Fig.~1 in~\cite{10.1063/5.0238408}, which has been published under the CC BY-NC license (see \url{https://creativecommons.org/licenses/by-nc/4.0/}).}\label{fig:schematic} 
    \end{figure}

    In this section we state the kinetic equation and the corresponding optimization problem that we investigate. The kinetic problem is given by the Vlasov-Poisson system with a source term that models the particle beam. The optimization problem is to find the best beam profile that minimizes the total electric energy in the plasma. A simple schematic of the setup we used is shown in figure~\ref{fig:schematic}.

    \subsection{The Kinetic Beam Heating Problem}
        The underlying kinetic equation is the Vlasov-Poisson system, which here describes the collisionless evolution of the electron distribution function $f$. We have
        \begin{equation}\label{eq:1}
            \partial_t f(t,x,v) + v \partial_x f(t,x,v) - E \partial_v f(t,x,v) = S(t,x,v) \qquad \text{ with }\qquad E = - \partial_x \phi.
        \end{equation}
        The potential $\phi$ is determined by solving $-\partial_{xx}\phi = \rho$, where $\rho = \frac{1}{L}\int f d(x,v) - \int f d(v)$ is the charge density. Note that $\int f d(x,v)$ models the positively charged stationary ion background. Since the source term $S$ adds particles to the plasma, this time-dependent term ensures that the plasma stays neutral overall. The initial condition for all our simulation is a homogeneous Maxwellian.
        The spatial domain is of size $L$, i.e. $x \in [0,L]$, with periodic boundary conditions and for the velocity we have $v \in [-6,6]$.

        The source term $S$ models the particle beam that enters the plasma. The beam is assumed to have a Gaussian cross section in space and is a (shifted) Maxwellian in velocity. More specifically, we have 
        \begin{align}\label{eq:2}
            S(t,x,v) = I(t) \frac{1}{\sqrt{2\pi}\sigma} \exp\left[ \frac{-(x-x_0)^2}{2\sigma^2}\right] \frac{1}{\sqrt{2\pi}v_{thb}} \exp\left[ \frac{-(v-v_b)^2}{2v_{thb}^2}\right],
        \end{align}
        where $x_0$ is the center of the beam. However, since the spatial domain is periodic the location of $x_0$ is irrelevant. In the numerical simulations we use $x_0=\frac{L}{2}$. The spatial extent of the beam is controlled by $\sigma$ and $v_{thb}$ is the thermal velocity of the beam. The average beam velocity is denoted by $v_b$.

        Injecting a high-speed beam into a plasma modifies the velocity distribution function by introducing an additional peak, or \enquote{bump}, at higher velocities. It is well established that the presence of this secondary peak can introduce instabilities in the plasma. This effect is commonly known as the bump-on-tail instability and has been extensively studied, see, e.g.~\cite{Denavit1985,Yamagiwa1987,Carlevaro2015,Zeeland2021}. 

        Within linear theory, the bump-on-tail instability is most commonly analyzed by assuming an initial distribution that is a linear combination of two Maxwellians. In this case, an instability develops as long as there is a positive slope of the distribution function at velocity $\omega / k$, where $\omega$ and $k$ are the frequency and wave-vector, respectively, of a Langmuir wave. Since Langmuir waves (or plasma oscillations) have a dispersion relation where $\omega \approx 1$ this means we require that our system supports $k = m 2\pi/L$, where $m\in \mathbb{N}$, such that $\omega / k$ falls into the velocity range where the slope is positive. In this case the plasma can transfer energy to the electric field which results in an exponential increase of the electric energy. In the linear theory this growth continues forever. However, in numerical simulations or experiments nonlinear effects eventually take over and lead to a saturation of the instability~\cite{D_G_Swanson_2003}.

        More realistically, the beam is continuously injected into the plasma over some time. Our model takes this into account via the source term. This situation is difficult to analyze analytically. However, for a uniform intensity beam, i.e. $I(t)=\text{const}$, the behaviour of the instability is very similar; see, the left of Figure~\ref{fig:bi_ee_plasma1}. As has been pointed out in~\cite{10.1063/5.0238408}, this also gives us the possibility to consider non-uniform beam profiles, i.e.\ different $I(t)$, in order to reduce or suppress the instability.

    \subsection{The Optimization Problem}

        The objective of this paper is to identify beam profiles that reduce or suppress the instability, while simultaneously preserving the total amount of energy injected in the plasma. 

        In principle, the optimizer could try to find an arbitrary function $I(t)$ that achieves this goal. However, this is not desirable as there are usually experimental constraints on the beam profile. In addition, it makes the optimization problem extremely high-dimensional and thus hard to solve. Therefore, $I(t)$ is commonly parametrized using a smaller number of parameters.

        This can be done in multiple ways, in~\cite{10.1063/5.0238408} the beam profile is parametrized by a combination of sines and cosines with different frequencies, i.e.\ we choose
        \begin{align}\label{eq:3}
            I(t) = \max \left( I_0 + \sum_{k=1}^{K-1} a_k \cos(k\omega t) + \sum_{k = 1}^{K} b_k \sin(k \omega t), 0 \right)
        \end{align}
        for $t \in [0,t_b]$, where $t_b$ is the beam duration, and $I(t)=0$ otherwise. The parameters $a_k, b_k$ and $\omega$ are varied through the optimization process. We denote the average beam intensity by $I_0 = 0.05$. We do not vary any of the other parameters from~\eqref{eq:2}. This leaves us with $2K$ parameters in total.

        There are two advantages of this approach. First, it produces a smooth beam profile that avoids rapid change. Second, it yields a function that is periodic in time. The latter, in particular, means that we do not need to pose the optimization problem for the entire duration the beam is switched ON, i.e.~in the time interval $[0,t_b]$, but can divide the problem into functions on periods of size $2\pi / \omega$. The period is also chosen by the optimization algorithm.

        It is important to note that the total mass injected by the beam throughout the simulation is always the same, i.e.
        \begin{equation}\label{eq:4}
            \int_0^{t_{final}} I(t) dt = t_b I_0.
        \end{equation}
        For the parametrization in equation~\eqref{eq:3} this is not exactly true as $I(t)$ is clipped to avoid negative values. However, if this occurs then the overall energy injected is larger which makes the instability more severe. Thus, the behavior is naturally penalized and thus not of much concern.

        The total electric energy is a metric for how severe the instability in the plasma is. If there is no instability, we simply get $E = 0$. Therefore, we choose to minimize the total time-integrated electric energy in the plasma as our Figure of Merit (FOM). More specifically, we use
        \begin{align}\label{eq:5}
            \text{FOM } = \frac{1}{2} \int_0^{t_{\text{final}}}\int_0^L E^2(t,x) dx dt,
        \end{align}
        where $t_{\text{final}}$ is the final time of the simulation. In principle, alternative metrics are possible, although they result in a somewhat different optimization and parameter-landscape. For a more thorough discussion we refer to~\cite{fom}.
            
        In the simulation we always normalize the results by the reference value. The reference value is the FOM~\eqref{eq:5} when the beam is switched ON with a constant intensity $I_0$ until time $t=t_b$, we call this the uniform beam. Normalizing allows us to more easily compare different beam profiles. Every value below $1$ means that the instability is reduced and in cases where the FOM is above $1$ this means that the instability is made more severe.

\begin{figure}[H]
            \centering
            \includegraphics[width = 1.0\linewidth]{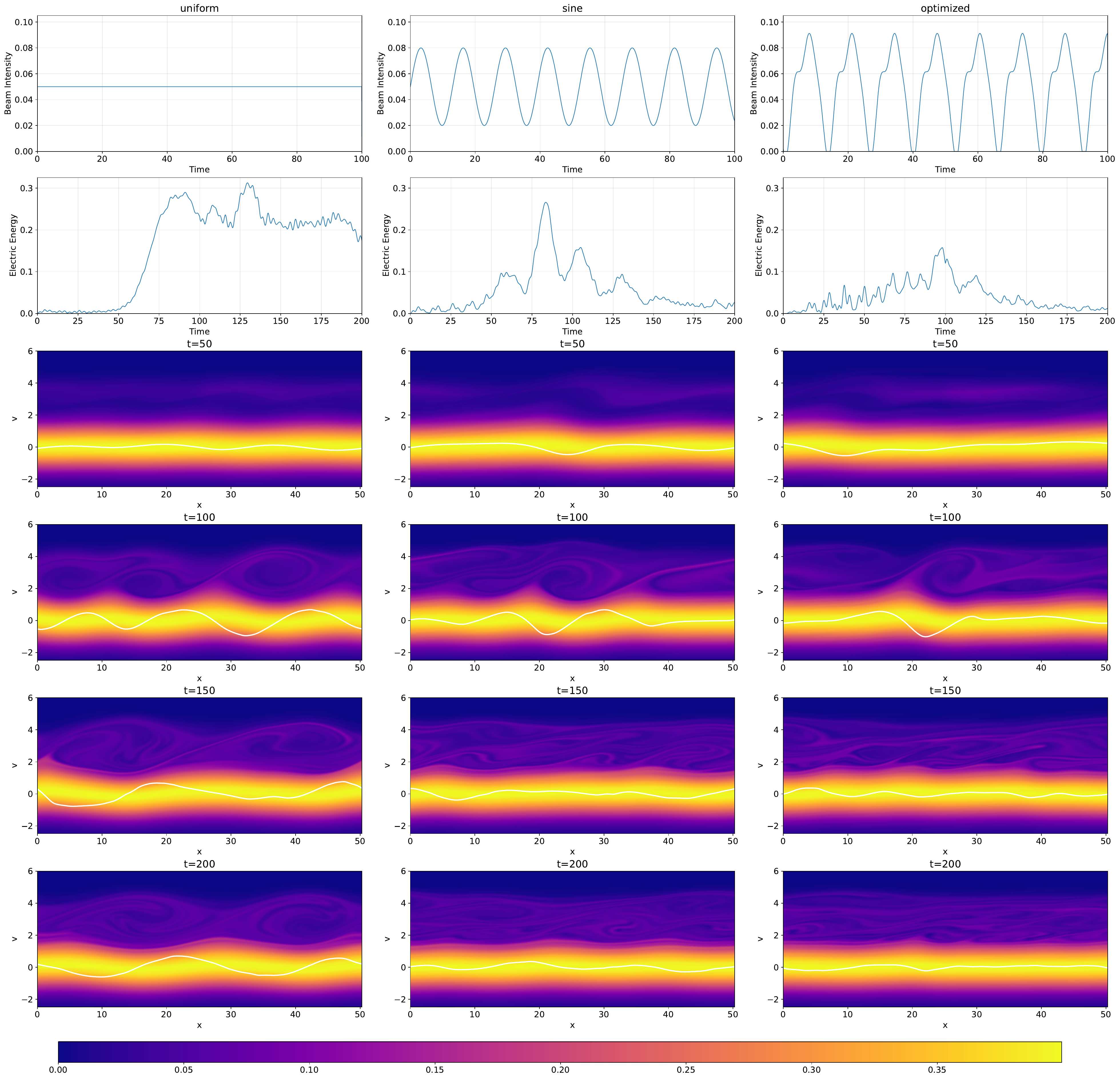}
            \caption{Comparison of beam profiles, electric energies and plasma distributions for times $t=50,100,150,200$. From left to right we compare the uniform beam ($FOM=1$), a sinusoidal beam ($FOM=0.42$) and the optimized beam found in~\cite{10.1063/5.0238408} ($FOM=0.28$). All parameters are listed in table~\ref{tab:2}.}\label{fig:bi_ee_plasma1}
        \end{figure}

        Figure~\ref{fig:bi_ee_plasma1} shows numerical simulations with different beam profiles. The beam profiles are shown in the top row  and correspond to the uniform beam (left), a sinusoidal beam (middle), and a modulated beam with $K=3$ (right). The modulated case corresponds to the parameters found in~\cite{10.1063/5.0238408} by an optimization algorithm. In the second row the electric energy as a function of time is shown. Below the electric energies we have $4$ snapshots of the plasma distribution function. We can clearly see that both the sinusoidal beam and the $K=3$ beam significantly reduce the FOM and thus the severity of the instability.

\begin{figure}[H]
            \centering
            \includegraphics[width = 1.0\textwidth]{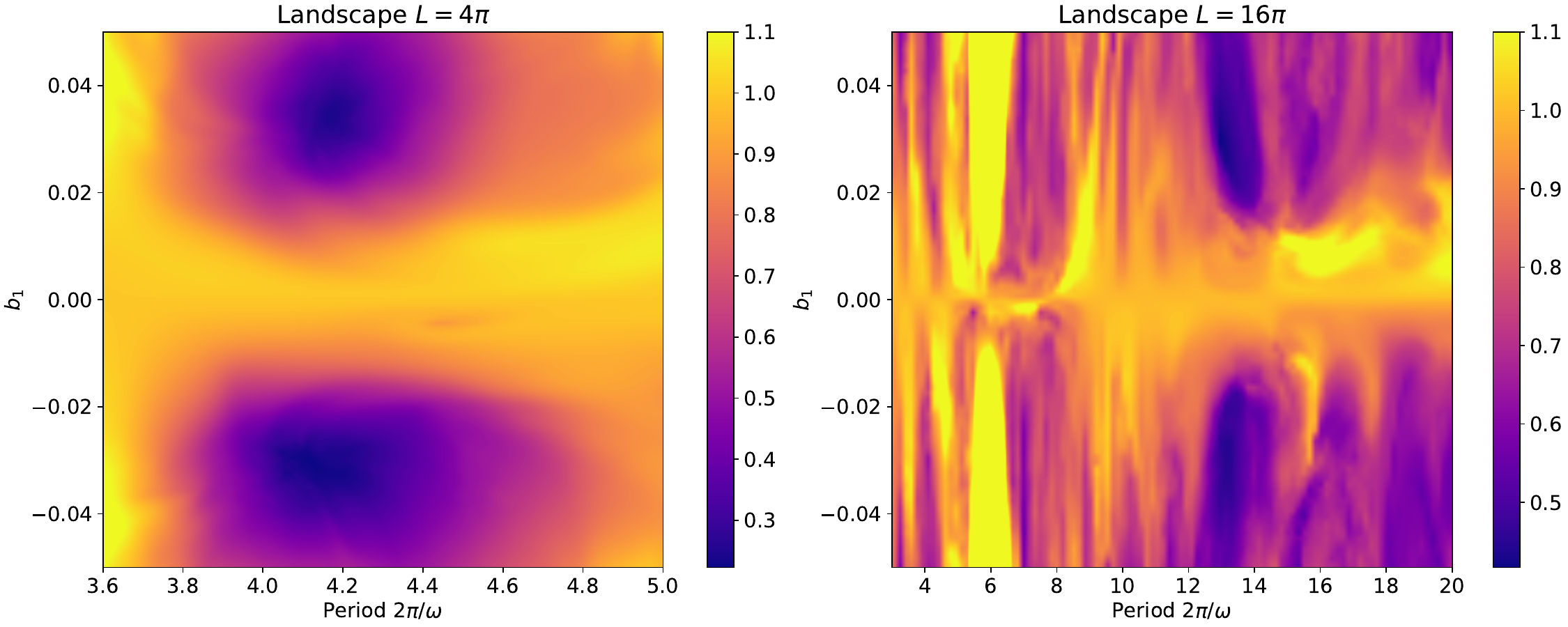}
            \caption{2D optimization landscapes for spatial domains $L = 4\pi$ (left) and $ L = 16\pi$ (right). The beam intensity is $I(t)=I_0 + {b_1} \sin(\omega t)$. We plot the FOM as a function of amplitude $b_1$ and the period $2\pi/\omega$.}\label{fig:L1_L2_landscape}
        \end{figure}

        The primary challenge in solving the optimization problem is that the underlying parameter landscape is very rugged with many local minima. To illustrate this we show a landscape plot for $L= 4\pi$ (LHS) and $L = 16\pi$ (RHS) in figure~\ref{fig:L1_L2_landscape} for the sinusoidal beam profile, i.e.~we consider $K = 1$ in~\eqref{eq:3}, which gives
        \begin{equation}
            I(t) = I_0 + b_1 \sin(\omega t).
        \end{equation}
        Note that the domain becomes visibly more rugged as the domain length becomes large compared to the Debye length (which in our unit system is $1$).

        In the following section, we introduce the optimization algorithms under consideration and evaluate their performance for the modulated case ($K=3$), i.e.~a six-dimensional optimization problem. Finally, we analyze the factors contributing to the success or failure of the respective algorithms.

\section{Comparison of the Optimization Algorithms}\label{sec:comparison}
    In this section we evaluate the performance of both local and global optimization algorithms for the problem described in section~\ref{sec:problem}. To do this we first introduce the considered algorithms: L-BFGS-B, Basinhopping, Bayesian Optimization (BO), CMA-ES, Particle Swarm Optimization (PSO) and Differential Evolution (DE). We then evaluate their performance and explain the reason for their failure or success.

    \subsection{Description of the Algorithms}
        \textbf{Limited-memory BFGS with Box constraints (L-BFGS-B)} is a quasi-Newton and hence local optimization algorithm designed for solving large-scale optimization problems with bound constraints~\cite{L-BFGS-B_paper}. In our numerical experiments we used the implementation from SciPy \texttt{(scipy.optimize.minimize)}~\cite{scipy,lbfgsb_implementation,lbfgsb_implementation1}. The L-BFGS-B algorithm iteratively performs the following steps:
        \begin{enumerate}
            \item \textbf{Gradient Evaluation}: Compute the gradient of the objective function at the current point.
            \item \textbf{Hessian Approximation Update}: Update the approximation of the inverse Hessian matrix using the limited-memory approach.
            \item \textbf{Search Direction Calculation}: Calculate the search direction based on the gradient and the Hessian approximation.
            \item \textbf{Line Search}: Perform a line search along the search direction to find an optimal step size.
            \item \textbf{Update Solution}: Update the current solution based on the search direction and step size.
            \item \textbf{Iteration}: Repeat the process until convergence or a stopping criterion is met.
        \end{enumerate}
        
        \textbf{Basinhopping} is a global optimization algorithm that combines local optimization with random perturbations to escape local minima~\cite{basinhopping_paper}. For the numerical experiments we used the algorithm from SciPy \texttt{(scipy.optimize.basinhopping)}~\cite{scipy,basinhopping_implementation}.~The algorithm iteratively performs the following steps:
        \begin{enumerate}
            \item \textbf{Random Perturbation}: A random perturbation is applied to the current solution to explore the search space and potentially move to a different basin of attraction.
            \item \textbf{Local Optimization}: Starting from the perturbed point, a local optimization method (e.g.~L-BFGS-B) is applied to find a local minimum of the objective function.
            \item \textbf{Acceptance Criterion}: The newly obtained solution is evaluated. If it has a lower objective function value than the old solution, it is accepted as the new current solution. Otherwise, it may still be accepted with a certain probability, to allow for exploration.
            \item \textbf{Iteration}: The process is repeated until a stopping criterion is met, such as a maximum number of iterations or convergence.
        \end{enumerate}

        \textbf{Bayesian Optimization (BO)} is a technique for optimizing expensive objective functions~\cite{bayes}. It is particularly useful when the objective function is costly to evaluate, as it aims to find the global minimum with as few evaluations as possible. The core idea behind BO is to build a probabilistic model of the objective function and use this model to make informed decisions about where to sample next. For the implementation of BO we used botorch and gpytorch packages~\cite{botorch,gpytorch}.The process of BO involves the following steps:
        \begin{enumerate}
            \item \textbf{Surrogate Model Construction}: A surrogate model, often a Gaussian Process (GP), is constructed to approximate the parameter landscape based on the observed data points. The surrogate model provides a mean prediction and uncertainty estimate for the objective function at any given point in the search space.
            \item \textbf{Acquisition Function Optimization}: An acquisition function is defined to guide the search for the next sampling point. The acquisition function balances exploration (sampling in areas with high uncertainty) and exploitation (sampling in areas with low predicted values).
            \item \textbf{Sampling and Updating}: The next point to evaluate is determined by optimizing the acquisition function. The objective function is then evaluated at this point, and the surrogate model is updated with the new data.
            \item \textbf{Iteration}: The process is repeated until a stopping criterion is met, such as a maximum number of iterations or a convergence threshold.
        \end{enumerate}

        \textbf{CMA-ES (Covariance Matrix Adaptation Evolution Strategy)} is a stochastic optimization algorithm belonging to the family of evolutionary strategies~\cite{CMA-ES}. It is particularly effective for solving complex, non-linear, and non-convex optimization problems where gradient information is unavailable. CMA-ES adapts a multivariate normal search distribution to efficiently explore the parameter space and learn correlations between variables. For the implementation of the CMA-ES algorithm we used the package pycma~\cite{cma_package}. The algorithm iteratively performs the following steps:
        \begin{enumerate}
        \item \textbf{Sampling}: Generate a population of candidate solutions from a multivariate normal distribution.
        \item \textbf{Evaluation of the Objective Function}: Evaluate the objective function for each candidate solution.
        \item \textbf{Selection}: Select the best-performing candidates.
        \item \textbf{Distribution Update}: Update the mean of the search distribution based on the selected candidates.
        \item \textbf{Covariance Matrix Adaptation}: Adapt the covariance matrix to capture promising search directions and parameter correlations.
        \item \textbf{Step-Size Adaptation}: Adjust the global step size to balance exploration and exploitation.
        \item \textbf{Iteration}: Repeat the process until a stopping criterion is met.
        \end{enumerate}

        \textbf{Particle Swarm Optimization (PSO)} is a population based optimization algorithm inspired by the social behavior of birds and fish~\cite{PSO}. PSO does not need gradient information. In PSO, a swarm of particles (candidate solutions) moves through the search space, adjusting their positions based on their own experience and the experience of the best solution so far. For the numerical experiments with PSO we used the pyswarm package~\cite{pyswarm_package}.The PSO algorithm iteratively performs the following steps:
        \begin{enumerate}
            \item \textbf{Initialization}: Initialize a swarm of particles with random positions and velocities in the search space.
            \item \textbf{Evaluation of the Objective Function}: Evaluate the FOM for each particle based on its current position.
            \item \textbf{Update Personal Best}: For each particle, update its personal best position if the current position leads to a lower objective function value.
            \item \textbf{Update Global Best}: Identify the global best position among all particles in the swarm.
            \item \textbf{Velocity Update}: Update the velocity of each particle based on its personal best position, the global best position, and its current velocity.
            \item \textbf{Position Update}: Update the position of each particle based on its updated velocity.
            \item \textbf{Iteration}: Repeat the process until a stopping criterion is met.
        \end{enumerate}

        \textbf{Differential Evolution (DE)} is a population based evolutionary optimization algorithm that has proven its value in solving complex optimization problems~\cite{StornPrice}. Similar to CMA-ES and PSO, DE does not require gradient information and is therefore suitable for optimization problems with non-linear, non-convex, or non-differentiable objective functions. DE evolves a population of candidate solutions by generating new solutions through mutation and crossover operations. There are different strategies for the crossover in DE, in this work, we use the DE/best/1/bin strategy. For the implementation of DE we used SciPy \texttt{(scipy.optimize.differential\_evolution)}\cite{scipy}. The algorithm iteratively performs the following steps:
        \begin{enumerate}
            \item \textbf{Initialization}: Initialize a population of individuals randomly in the search space and evaluate their objective function values.
            
            \item \textbf{Selection of Individuals}: Select the best individual, and for every individual, two further distinct individuals from the population.
            
            \item \textbf{Mutant Vector Generation}: Create a mutant vector using the best individual and the scaled difference between the two selected individuals. The preselected mutation factor controls the magnitude of this difference and therefore determines the size of the perturbation applied to the best individual.
            
            \item \textbf{Crossover}: Perform a crossover between the mutant vector and the current individual to create a trial vector. For every parameter, the corresponding value is selected either from the mutant vector or from the current individual, according to the chosen crossover probability. At least one parameter value is always taken from the mutant vector.
                        
            \item \textbf{Selection}: Compare the FOM of the trial vector with that of the current individual. If the trial vector has a better FOM, it replaces the current individual in the population. Otherwise, the current individual is retained.
            
            \item \textbf{Iteration}: Repeat the process until a stopping criterion is met.
        \end{enumerate}

        \textbf{Enhanced Differential Evolution}
        To investigate whether the performance of DE can be improved by adding refinement and initialization strategies, three enhanced variants of the standard DE algorithm were implemented and analyzed.:
        \begin{itemize}
            \item \textbf{DE\_LG} Use a local polisher (L-BFGS-B) after every generation to refine the candidate solutions.
            \item \textbf{DE\_DE} DE is supplemented by periodic local DE searches. Every n-th generation, a small subset of individuals is selected and perturbed to form initial populations for separate local DE runs. Each selected individual is then replaced by the best solution obtained from the corresponding local search.
            \item \textbf{DE\_more\_init} DE is initialized from a larger set of randomly sampled individuals than required for the initial population. The individuals with the best objective-function values are then selected to form the initial population.
        \end{itemize}
        Among the three variants, only \textbf{DE\_more\_init} resulted in a measurable improvement over the standard DE algorithm. The performance of all the listed optimization approaches is evaluated in the following section.

    \subsection{Results\label{sec:smooth-beam-results}}
        In this section we present the results of the different optimization algorithms and then compare the performance of the algorithms by the FOM as a function of cost. The cost is the number of evaluations of the solver code. We used a solver developed in C++ to solve the Vlasov-Poisson system~\eqref{eq:1} that processes multiple beam profiles in batches. The implementation uses the Kokkos performance‑portability framework, enabling efficient execution on modern systems, including GPUs. We apply a splitting‑based semi‑Lagrangian method with ninth‑order Lagrange interpolation. This solver was already introduced and used in~\cite{10.1063/5.0238408}.
        For comparing the algorithms we chose to run the problem in a bounded six-dimensional parameter domain ($K=3$) with $L = 16\pi$. The used parameters are given by
        \begin{align}\label{eq:parameterset}
            t_f &= 200, t_b = 100, \Delta t = 0.02, L = 16\pi, N_x = 256, N_v = 256 \nonumber\\
            \sigma &= 0.8, v_b = 3.5, v_{thb} = 0.5, I_0 = 0.05, P = \frac{2\pi}{\omega} \in [12.5, 14.0], (a_k, b_k)_{k=1}^K \in [-0.05, 0.05]^{2K}
        \end{align}
        The parameters we vary in the optimization process are the period $P$ and the coefficients $a_1, a_2, b_1, b_2, b_3$. The final simulation time is denoted by $t_f$, while $t_b$ represents the duration for which the beam remains switched ON.~The time step size and the number of grid points in space and velocity are denoted by $\Delta t$, $N_x$ and $N_v$ respectively.
        For L-BFGS-B, we generated $1000$ initial points by Latin Hypercube Sampling and performed local optimization from each starting point. We used $\texttt{gtol}=10^{-3}$, $\texttt{ftol}=10^{-3}$, $\texttt{maxiter}=50$, and $\texttt{maxfun}=150$.

        Basinhopping was run as $60$ independent searches in parallel, each initialized with a uniformly random point within the bounds. The runs used three groups of $20$ with different step sizes. For local refinement we used the L-BFGS-B minimizer which was limited to $4$ iterations and $8$ function evaluations.
        
        For Bayesian optimization, we first evaluated $300$ uniformly sampled points and subsequently evaluated candidate points in batches of $50$. We used a variational Gaussian process and $q$-Log expected improvement as the acquisition function.

        For CMA-ES, the six-dimensional parameter domain was normalized, and the algorithm was initialized at the midpoint of this domain. We used an initial step size of $\sigma_0=0.2$ and a population size of $60$.

        Particle Swarm Optimization was initialized with $60$ particles. The particle positions were generated uniformly at random within the prescribed parameter bounds, and the optimization was run with acceleration coefficients $c_1=c_2=1.5$ and inertia weight $w=0.7$.

        Differential Evolution used a Latin Hypercube initial population of $60$ individuals. The mutation factor was set to $F=(0.5,1)$ and the crossover probability to $CR=0.7$. The DE/best/1/bin strategy was employed for generating trial vectors.
        
\begin{figure}[H]
            \centering
            \includegraphics[width = 1.0\textwidth]{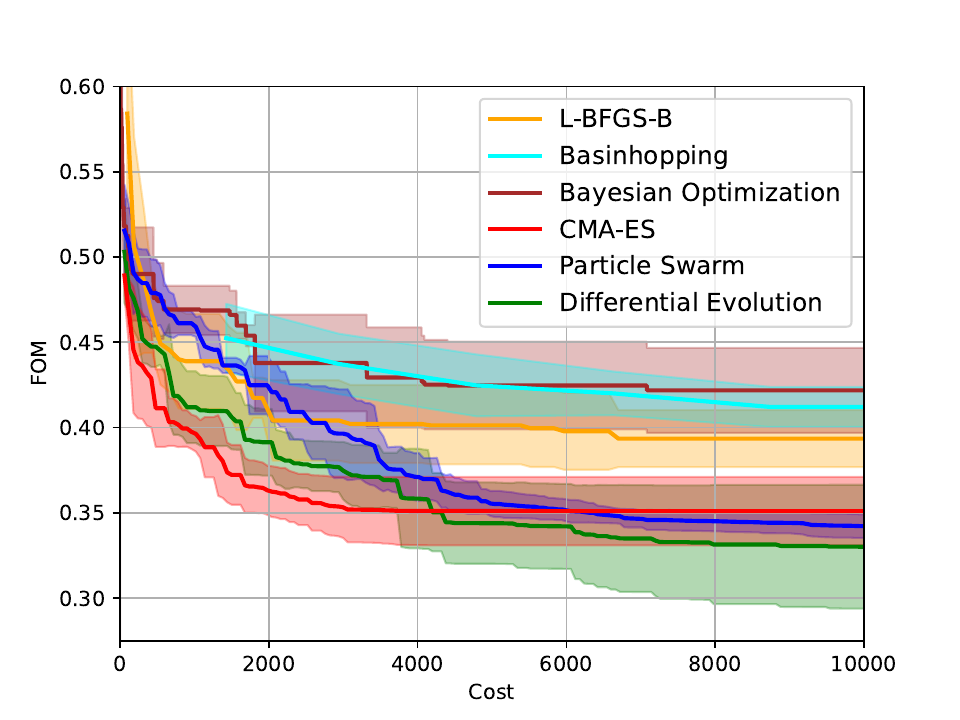}
            \caption{Comparison of all analyzed algorithms for the setup described in~\eqref{eq:parameterset}. We plot the mean and standard deviation as a function of cost. For each algorithm we ran $5$ repetitions with different initial values and seeds.}\label{fig:std_all}
        \end{figure}

        In figure~\ref{fig:std_all} we can clearly see that the population based methods CMA-ES, PSO and DE perform best for our problem. CMA-ES converges the fastest as we can see in the bottom left corner of figure~\ref{fig:std_all}. This is particularly useful if our computational budget is limited. However, CMA-ES does not find the best parameter set as we can see on the right of figure~\ref{fig:std_all}. More precisely, after $10000$ function evaluations DE reaches the lowest FOM of $0.284$ with an average of $0.33$. In comparison, CMA-ES reaches down to a minimum of $0.325$ with an average of $0.351$. Both PSO and DE find better solutions to our optimization problem if we can afford a large number of function evaluations. PSO converges slower than DE and is also not able to reach the same minimum. PSO reaches its lowest value of $0.329$ with an average of $0.342$ after $10000$ function evaluations. However, DE has by far the widest spread of the population based methods. PSO is the most consistent method. The least successful methods of the comparison are L-BFGS-B, Basinhopping and BO, since they do not get close to the minimum within the computational budget allotted. Note that the lines in figure~\ref{fig:std_all} start after finishing the first iteration of the corresponding optimization algorithm, which explains why they do not start at the same cost. In the following we will analyze these results in more detail.
        \begin{center} \textbf{Local Optimization fails} \end{center}
            The L-BFGS-B algorithm was set up with Latin Hypercube Sampling for the initial points in the parameter space. L-BFGS-B is very effective for smooth optimization problems, but it may struggle with rough landscapes that have many local minima. The algorithm relies on gradient information, which can be misleading in such landscapes, leading to convergence to local minima rather than the global minimum. Additionally, the line search may not be effective in navigating through rugged landscapes, further hindering its performance in this context. 
            \begin{figure}[H]
                \centering
                \includegraphics[width = 1.0\textwidth]{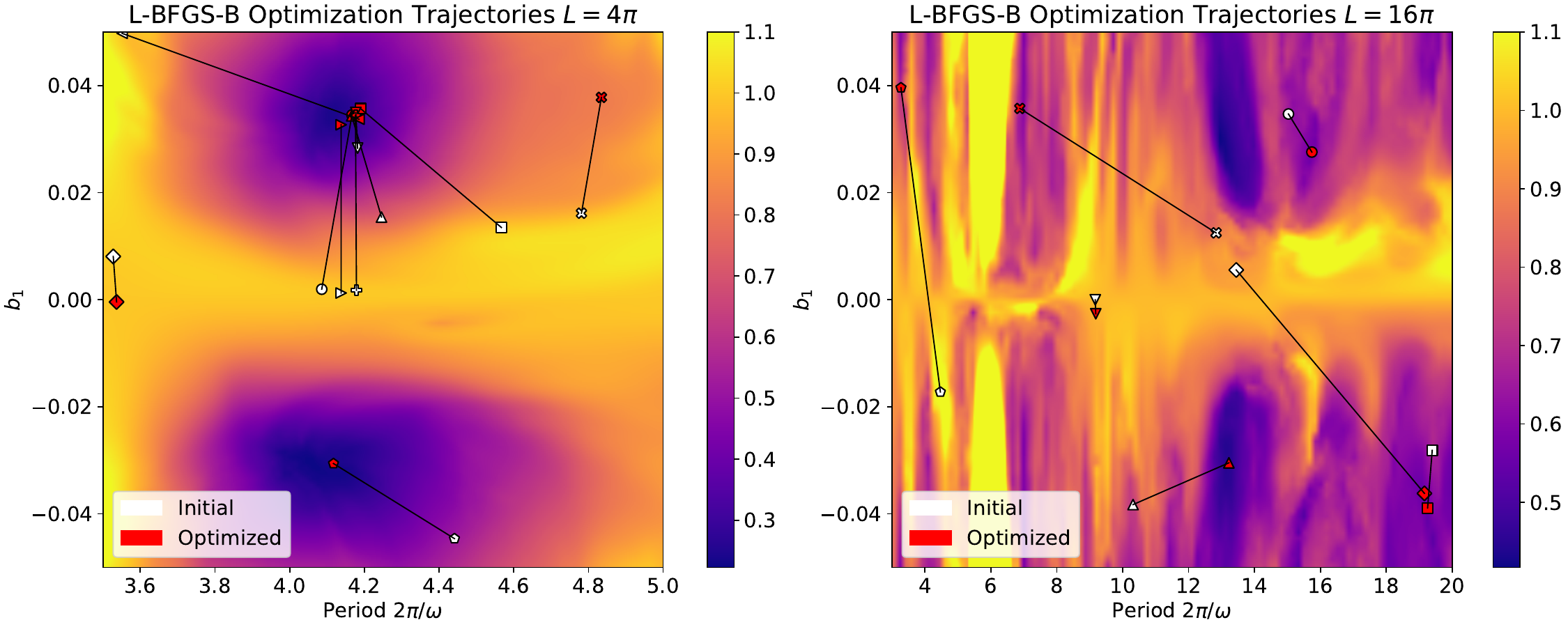}
                 \caption{The L-BFGS-B trajectories for spatial domains $L = 4\pi$ (LHS)  and $L = 16\pi$ (RHS) in 2D parameter space ($K = 1$) are shown. The different shapes represent different L-BFGS-B runs, where the white shapes are the initial parameters and the red shapes are the final parameters. The lines connect the initial and final points of each run.}\label{fig:L-BFGS-B_trajectories}
            \end{figure}
            We consider optimization with L-BFGS-B in more detail in figure~\ref{fig:L-BFGS-B_trajectories}. The L-BFGS-B algorithm works very well in the $L = 4 \pi$ case, almost regardless of the initial starting point the algorithm converges to one of the two acceptable basins. However, in the $L = 16 \pi$ case the algorithm is very sensitive to the initial starting point and often moves in the wrong direction, even further away from the global minimum. 
            \begin{figure}[H]
                \begin{minipage}{0.5\textwidth}
                    \centering
                    \includegraphics[width = 0.9\textwidth]{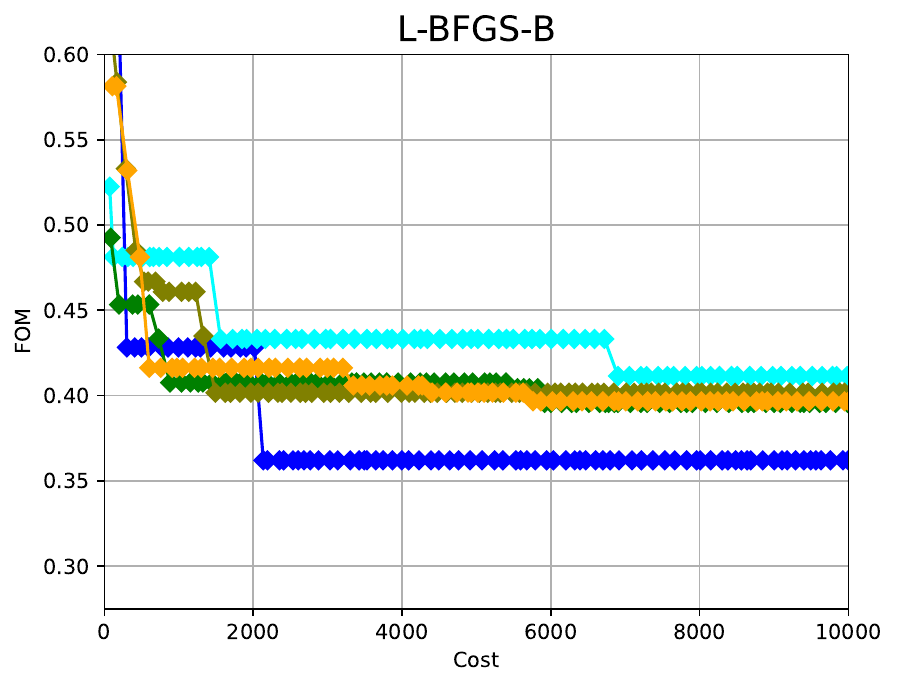}
                \end{minipage}\hfill
                \begin{minipage}{0.5\textwidth}
                    \centering
                    \includegraphics[width = 0.9\textwidth]{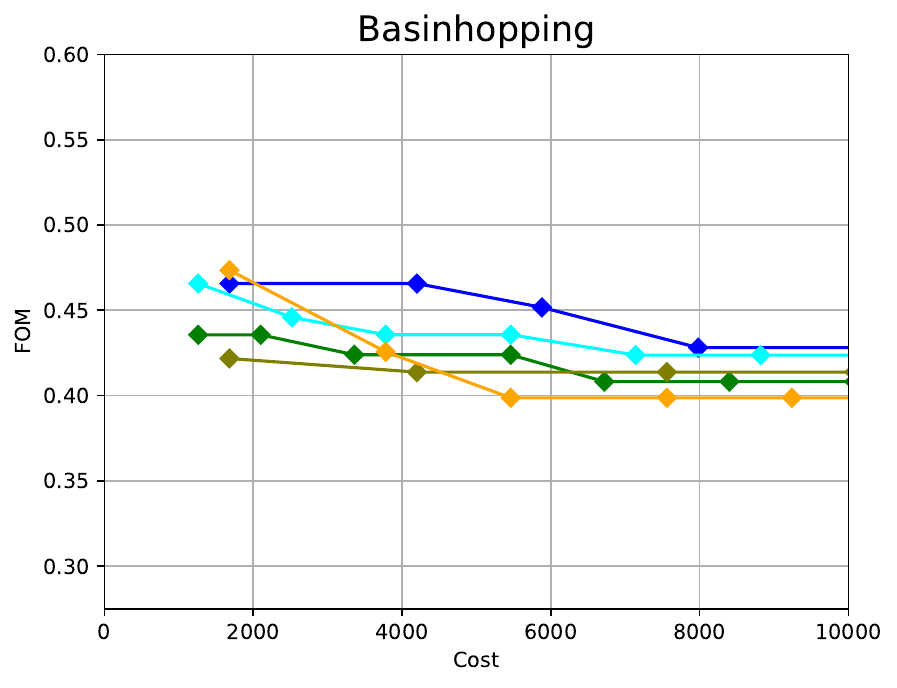}
                \end{minipage}
                \caption{(LHS) 5 L-BFGS-B runs with Latin Hypercube Sampling for the initial parameters. (RHS) 5 Basinhopping runs. For each Basinhopping run we started 60 algorithms independently in parallel with different initial sample points and iterated through generations. One line represents the best solution found so far of all 60 runs for each generation.}\label{fig:LBFGSB_basinhopping}
            \end{figure}
            Similar to L-BFGS-B, Basinhopping is not able to get close to the global minimum of the $L = 16 \pi$ landscape as we can see in figure~\ref{fig:LBFGSB_basinhopping}. Even though, in principle, Basinhopping is able to escape local minima this is not sufficient to get good performance for the problem studied here.
            
        \begin{center} \textbf{Surrogate methods scale badly with dimension} \end{center}

\begin{figure}[H]
                \begin{minipage}{0.5\textwidth}
                    \centering
                    \includegraphics[width = 0.9\textwidth]{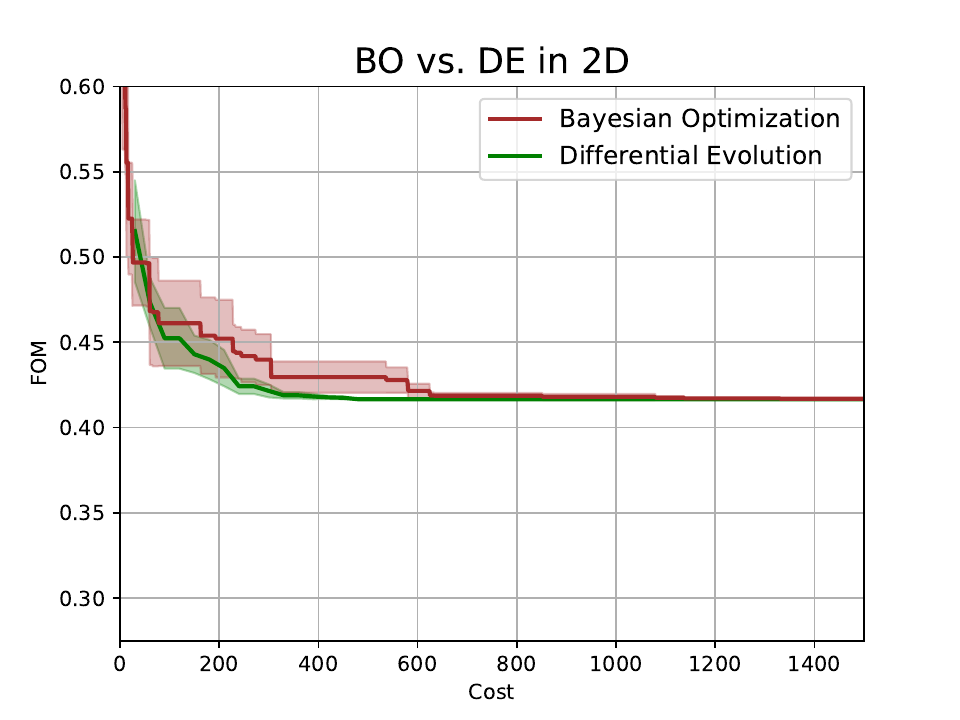}
                \end{minipage}\hfill
                \begin{minipage}{0.5\textwidth}
                    \centering
                    \includegraphics[width = 0.9\textwidth]{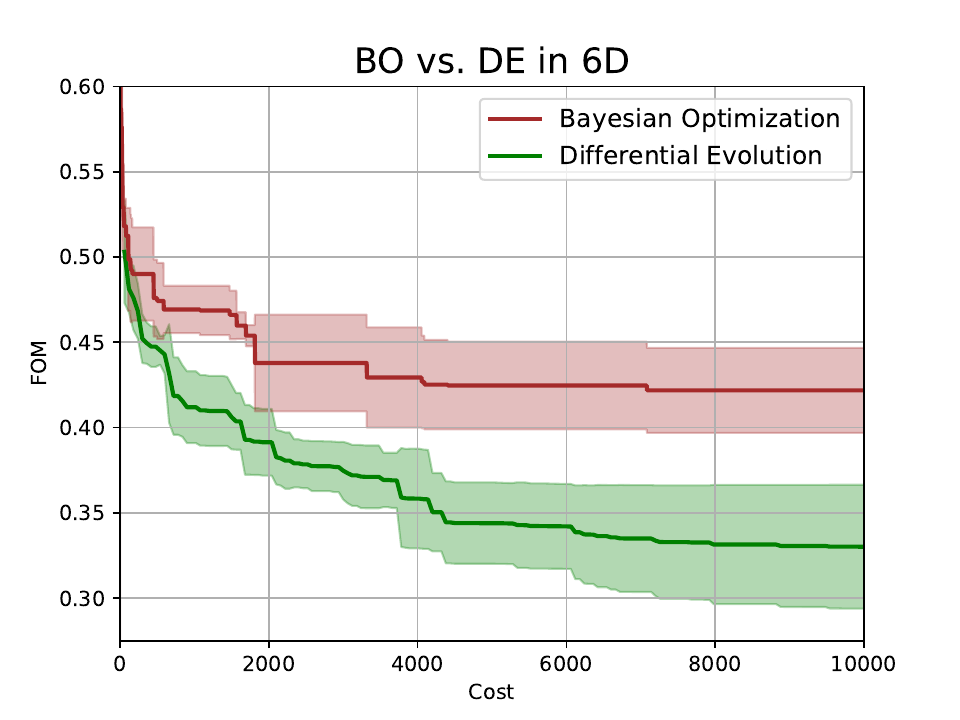}
                \end{minipage}
                \caption{BO does not scale well with dimensionality. We plot the mean and standard deviation of the algorithms as a function of cost. For the $2$D plot (LHS) we used $4$ simulations and for the $6$D plot (RHS) we used $5$ simulations. On the left-hand-side we see similar results for DE and BO.~However, on the right-hand-side we compare the performance in $6$D and clearly see the disadvantage of BO.}\label{fig:BO_vs_DE}
            \end{figure}

            Let us now compare Bayesian Optimization to Differential Evolution in two different cases. In both cases we consider $L = 16\pi$. In the first setting, we use $K = 1$, in~\eqref{eq:3}, while in the second setting, we use $K = 3$. Hence, we have two optimization problems with dimensions $2$ and $6$, respectively.
            The results are shown in figure~\ref{fig:BO_vs_DE}. We clearly see that Bayesian Optimization is worse than DE in both the 2D and 6D cases. However, in the left plot, representing the 2D case, we can see that BO is able to find the same solution as DE with a bit more computational effort. In the right plot, representing the 6D parameter space, BO is not able to come close to DE.~This suggests that BO scales badly with dimensionality; i.e.~the surrogate is not able to capture the roughness of the landscape in higher dimensions.

        \begin{center} \textbf{Population Based Methods} \end{center}

\begin{figure}[H]
                \centering
                \begin{tabular}{cc}
                    (A) & (B) \\
                    \includegraphics[width=0.47\textwidth]{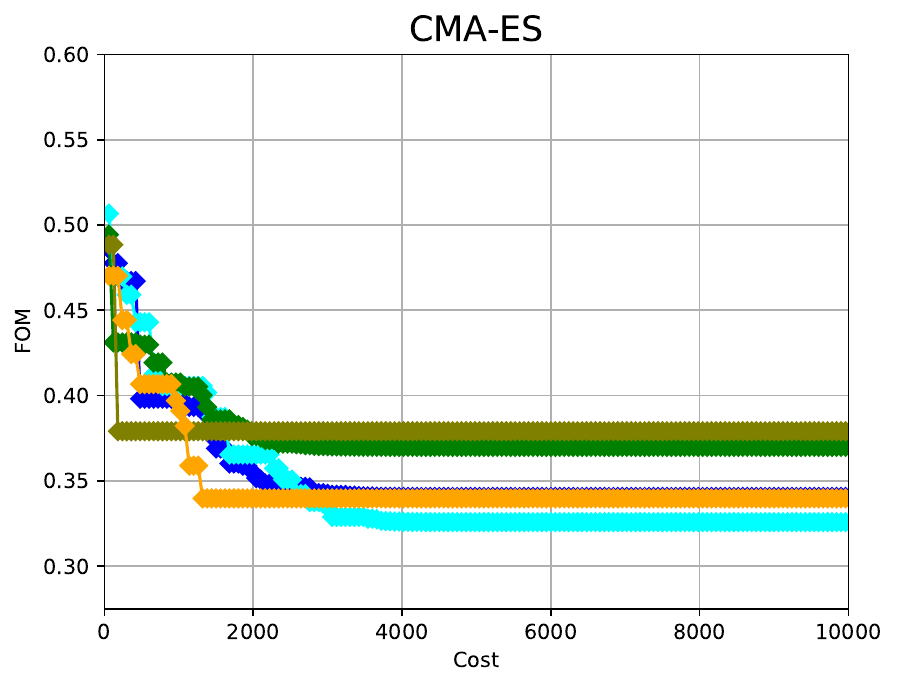} &
                    \includegraphics[width=0.47\textwidth]{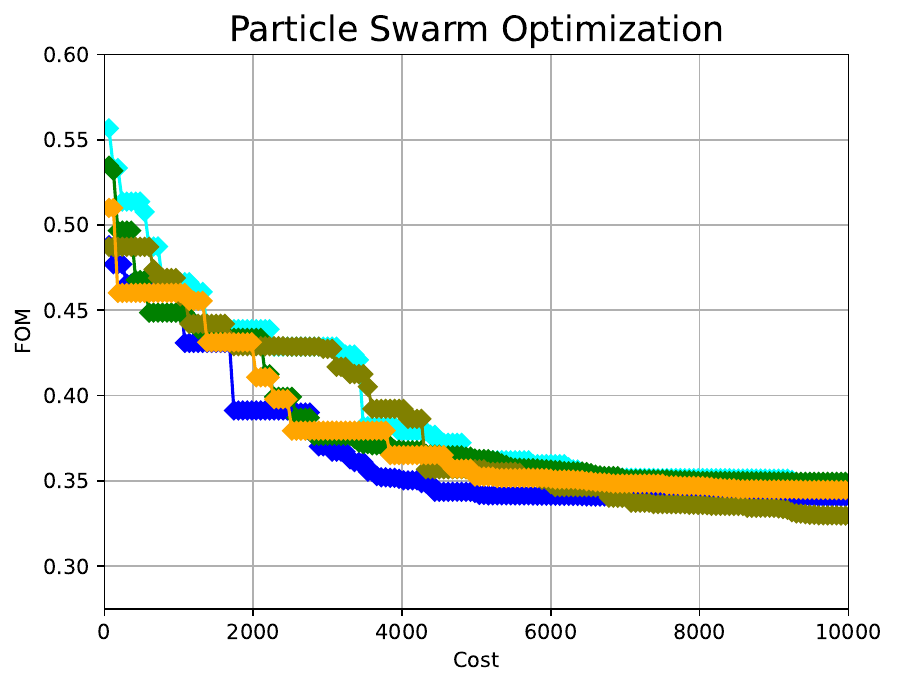} \\
                    (C) & (D) \\
                    \includegraphics[width=0.47\textwidth]{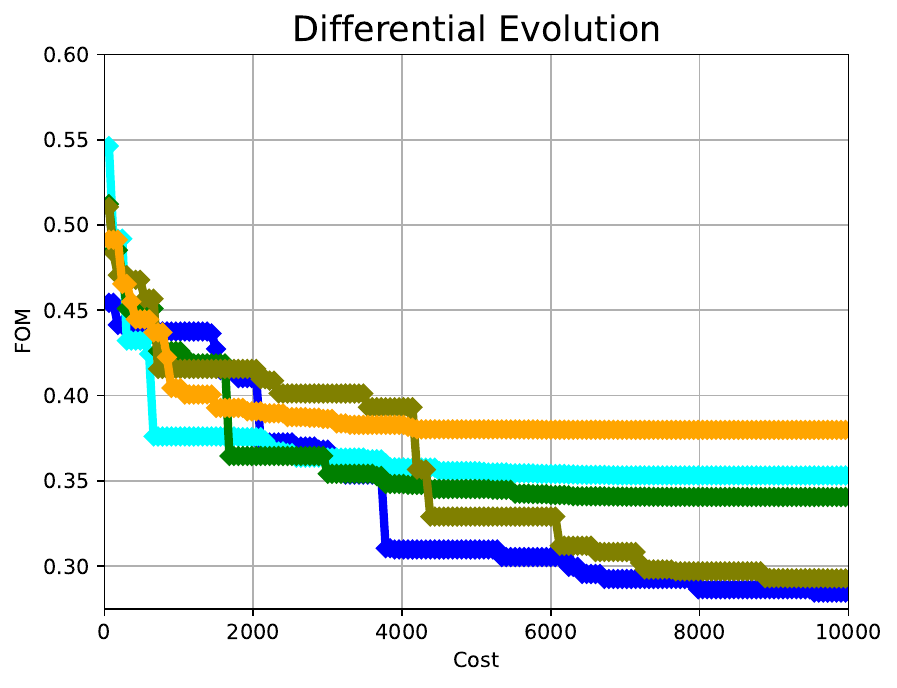} &
                    \includegraphics[width=0.47\textwidth]{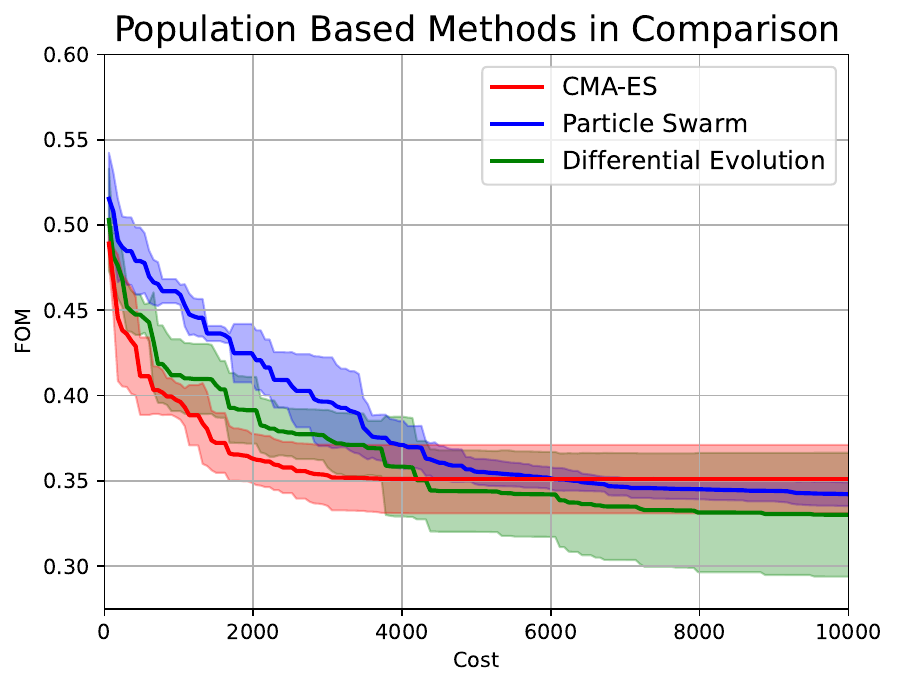}
                \end{tabular}
                \caption{
                    We analyze the performance of the population based methods, (A) CMA-ES, (B) Particle Swarm Optimization, (C) Differential Evolution, (D) Mean and standard deviation of the population based methods. For all three methods we compare $5$ runs with different initial populations and seeds.}\label{fig:population_based}
            \end{figure}

            The very clear winners of the comparison are the population based methods, e.g.~CMA-ES, Particle Swarm Optimization and Differential Evolution. We will now investigate their performance in more detail. We can clearly see from plot (D) in figure~\ref{fig:population_based} that DE has the widest spread. This implies that the algorithm is less predictable and can converge to worse but also to better solutions than the other population based algorithms. This can also be seen in plot (C), where we compare five different DE runs with each other. In these runs only the initial population is changed. In plot (C) the best performance of DE reaches a value of $0.285$ for our FOM after $10000$ cost, which is significantly better than both CMA-ES and Particle Swarm Optimization that attain $0.325$ and $0.329$, respectively. The large spread indicates that in order to obtain reliable results with DE, we might have to run the algorithm multiple times, which is expensive. In plot (B) we can see that PSO is very consistent and converges to a good solution for each of the five runs. CMA-ES initially converges very fast, after $3000$ cost there is no more significant change in the suggested minimum as can be seen in plot (A). CMA-ES is thus a good choice for optimization with a lower computational budget. Its consistency lies in between DE and PSO.

        \begin{center} \textbf{Enhanced Differential Evolution} \end{center}
            In this section we will now explain the enhancements of DE mentioned in section~\ref{sec:comparison} and analyze their performance.
            
            For the DE\_LG scheme we used L-BFGS-B every generation on $3$ random individuals. This would allow the algorithm to make use of gradient information. We stopped L-BFGS-B if the improvement in each gradient step is less than $1\%$.
            
            For the DE\_DE scheme we take $2$ random individuals and the best candidate every fifth generation out of the population and start three DE runs of $10$ individuals each, lasting for $10$ generations. The first individual is initialized at the selected starting point $x_0$, while the remaining $9$ individuals are generated by independently perturbing $x_0$ in each dimension. Specifically, for each variable $j$, a perturbation of the form
                \begin{equation}
                    \delta_j = s_j u_j (b_j^{\max} - b_j^{\min})
                \end{equation}
            is added to $x_0$. Where $s_j \in \{-1,1\}$ is a uniformly sampled random sign, $u_j \sim \mathcal{U}(0.1,0.2)$ and $b_j^{\max}, b_j^{\min}$ are the upper and lower bounds for variable $j$. Thus, each coordinate is perturbed by a random amount between $10\%$ and $20\%$ of its feasible range, with equal probability of a positive or negative displacement. The best value of these sub-DE runs replaced the original candidates. This allows us to explore the local parameter space more thoroughly and end up with better candidates for the next generation.
            
            For the DE\_more\_init scheme we sampled $600$ individuals by Latin Hyper Cube sampling and chose the best $60$ as the initial population. 
            
\begin{figure}[H]
                \centering
                \includegraphics[width = 0.75\textwidth]{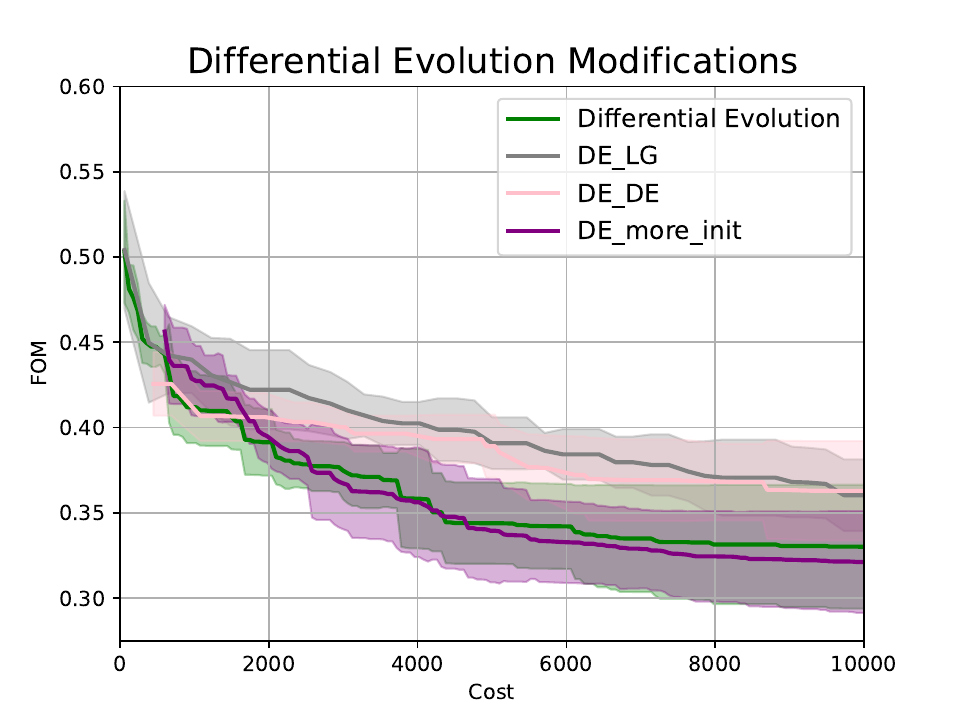}
                \caption{In this figure we compare the performance of the DE enhancements. For each scheme we used $5$ runs with different seeds and initial populations. The methods start out of sync since the initialization processes differ from each other.}\label{fig:DE_enhancements}
            \end{figure}

            As one can see in figure~\ref{fig:DE_enhancements} the only successful improvement of DE was made by sampling a larger initial population. The other two enhancements did not improve the performance of DE.~This is a bit surprising since in many problems local optimizers can increase the convergence rate. However, as we have seen gradient information often points the optimizer in the wrong direction, which harmed the performance of DE in our case. The DE\_DE scheme was also not successful, since the local DE runs oftentimes could not find better solutions than the initial points in the given budget, and hence did not improve the performance of the main DE run. Both of those point to the fact that local information is only of limited value for our problem.

            Let us note that DE\_more\_init does not give lower FOM values, but we reached the lowest FOM value faster and more consistently as compared to the original DE scheme. This is due to the fact that we have a better initial population and hence a better starting point for the optimization process.
            In the next section, we consider a different and more experimentally feasible parametrization of the beam profile, which we refer to as a Binary Beam. The optimization methods that proved successful for the previous parametrization are then applied to the Binary Beams, and the resulting performance is analyzed.

\section{Binary Beam}
\begin{figure}[H]
        \centering
        \begin{minipage}{0.49\textwidth}
            \centering
            \begin{tikzpicture}
            \node[inner sep=0] (beamplot) {\includegraphics[width=1\textwidth]{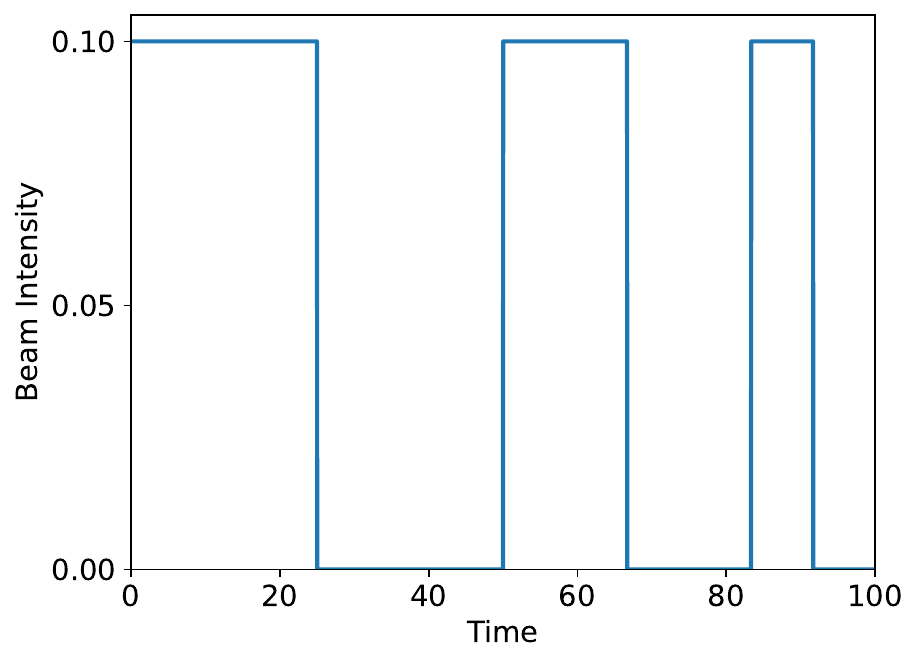}};
            \coordinate (h1start) at ($(beamplot.south west)+(11mm,-2mm)$);
            \coordinate (h1end) at ($(beamplot.south west)+(43.5mm,-2mm)$);
            \coordinate (h2start) at ($(beamplot.south west)+(43.5mm,-2mm)$);
            \coordinate (h2end) at ($(beamplot.south west)+(64mm,-2mm)$);
            \coordinate (h3start) at ($(beamplot.south west)+(64mm,-2mm)$);
            \coordinate (h3end) at ($(beamplot.south west)+(75mm, -2mm)$);
            \draw[decorate,decoration={brace,mirror,amplitude=4pt}] (h1start) -- (h1end)
               node[midway,below=6pt] {$H_1$};
            \draw[decorate,decoration={brace,mirror,amplitude=4pt}] (h2start) -- (h2end)
                node[midway,below=6pt] {$H_2$};
            \draw[decorate,decoration={brace,mirror,amplitude=4pt}] (h3start) -- (h3end)
                node[midway,below=6pt] {$H_3$};
        \end{tikzpicture}
        \end{minipage}
        \begin{minipage}{0.49\textwidth}
            \centering
            \includegraphics[width = 1\textwidth]{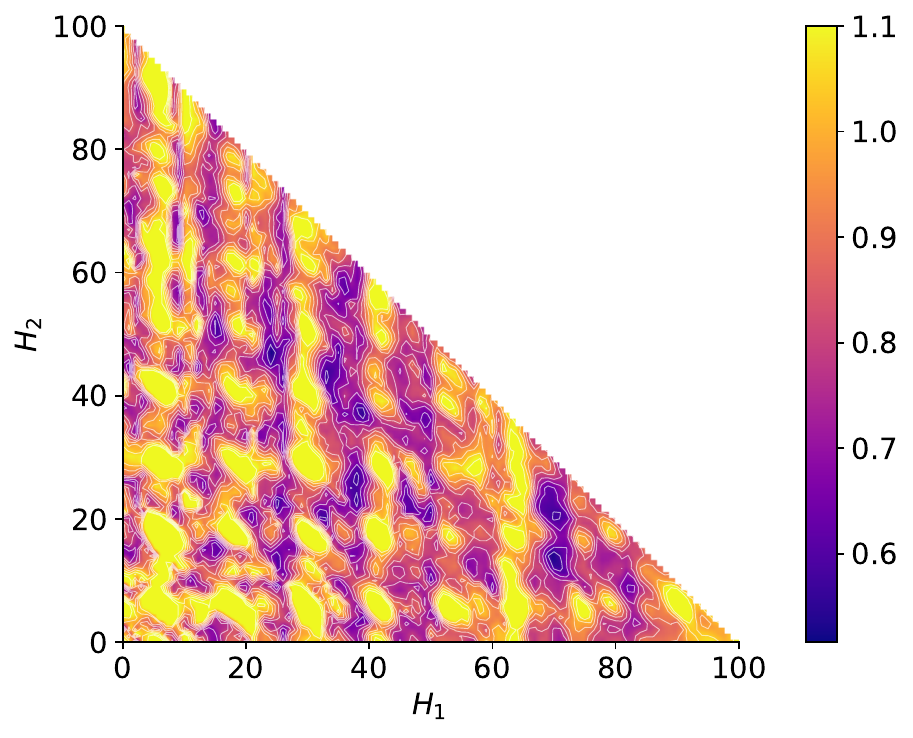}
        \end{minipage}\hfill
        \caption{(LHS) Example of a beam profile with parameters $t_1 = 0.75$, $t_2 = 0.5$, $t_3 = 0.25$. (RHS) The corresponding parameter landscape of the binary beam profile for $N=3$. The $x$ and $y$ axes represent $H_1$ and $H_2$ respectively. The third parameter $H_3$ is then uniquely determined. Note that since $H_1 + H_2 + H_3 = t_b$ has to be fulfilled, the upper triangle of the right plot corresponds to invalid parameter combinations.}\label{fig:binary_beam}
    \end{figure}

    In this section we consider the more restricted (but experimentally easier to realize) case of beam profiles that can only be ON (beam present) and OFF (no beam present). We will continue using all the same numerical parameters from~\eqref{eq:parameterset} and the $L=16\pi$ case. However, the binary beams will be parametrized by
    \begin{equation}
        \begin{bmatrix}
        t_1 & t_2 & t_3 & \cdots & t_{N}\\
        \end{bmatrix}
        \in [\epsilon, 1]^N.
    \end{equation}
    The $t_i$ give the relative length of each ON-OFF cycle. The duration of each cycle is then given by
    \begin{align}
        H_i =  \frac{t_b t_i}{\sum_{j = 1}^{N}t_j}.
    \end{align}
    In each cycle the beam is set to ON (i.e. $I(t)=2 I_0$) and set to OFF (i.e. $I(t)=0$) for a time interval of size $H_i/2$. Thus, the total energy of the beam that enters the plasma is still the same as for the uniform beam, as desired. The left of Figure~\ref{fig:binary_beam} shows an example of a beam profile with parameters $t_1 = 0.75, t_2 = 0.5, t_3 = 0.25$. Not all parameters in the parametrization are independent degrees of freedom as for the duration of each cycle $H_i$ only the relative size of the $t_i$ matters. However, this approach has the advantage that no constraints have to be added to the optimization problem.

    Note that we now optimize for the entire beam duration and not consider a periodic beam profile as we did before. Although, in principle, this would be possible for the binary beam case as well, we have observed worse performance following that approach. The downside of the non-periodic approach is that we now generally have a higher-dimensional optimization problem. However, we should note that since the individual cycles are very simple, the number of parameters is still manageable.

    Despite these differences, the optimization problem poses similar challenges. From the right-hand side of figure~\ref{fig:binary_beam} we can see that this parametrization leads, similar to the smooth beam profile, to a rugged optimization landscape.
    In this parametrization the number of cycles is fixed to $N$. However, it is unclear if a certain number of cycles is optimal. Thus, we will solve the optimization problem for different values of $N$ in the following.

    It is also desirable, from a practical point of view, to not allow for cycles with extremely short duration, as this might be challenging to realize experimentally. We can achieve this by setting a lower bound on the parameters, thus we look for $t_i \in [\epsilon ,1]$, for feasibly large values of $\epsilon$. This condition immediately leads to an upper and lower bound on the duration of the cycles $H_i$ given by
    \begin{align}
        H_i & \leq \frac{t_b}{(N-1)\epsilon + 1} , &&\text{ with equality if } (t_k)_{k=1}^{N} = ( \epsilon, \dots, 1, \dots, \epsilon )
    \end{align}
    and
    \begin{align}
        H_i &\geq \frac{t_b\epsilon}{N-1+\epsilon} {\approx \frac{t_b \epsilon}{N}}, &&\text{ with equality if } (t_k)_{k=1}^{N} = ( 1, \dots, \epsilon, \dots, 1 )\nonumber.
    \end{align}
    In the following we distinguish between results with $\epsilon = \frac{1}{50} = 0.02$ and $\epsilon = \frac{4}{25} = 0.16$. In table~\ref{tab:1} we summarize the upper and lower bounds for the duration of the cycles $H_i$, depending on the number of cycles $N$ and the lower bounds $\epsilon$.
    \setlength\extrarowheight{5pt}
    \begin{table}[H]
        \centering
        \begin{tabular}{c|cccc}
            $\epsilon$ $\backslash$ N & 5 & 10 & 15 & 20\\[1ex]
            \hline
            $\frac{1}{50}$ & [0.50, 92.59] & [0.22, 84.75] & [0.14, 78.13] & [0.11, 72.46]\\[1ex]
            $\frac{4}{25}$ & [3.84, 60.98] & [1.75, 40.98] & [1.13, 30.86] & [0.84, 24.75]\\[1ex]
        \end{tabular}
        \caption{Theoretical upper and lower bounds for the duration of the cycles $H_i$ for different values of $N$ and $\epsilon$. We have chosen $t_b=100$.}\label{tab:1}    
    \end{table}

    To distinguish between the different lower bounds we use different line styles in the plots, the results with full lines are computed with $\epsilon = \frac{1}{50}$ and the dashed lines represent the results with $\epsilon = \frac{4}{25}$. We call the latter, the smaller bounds case. We analyze these two cases in order to understand if having the ability to use shorter cycles is beneficial for minimizing the FOM.
    
\begin{figure}[H]
        \centering
        \includegraphics[width = 0.9\textwidth]{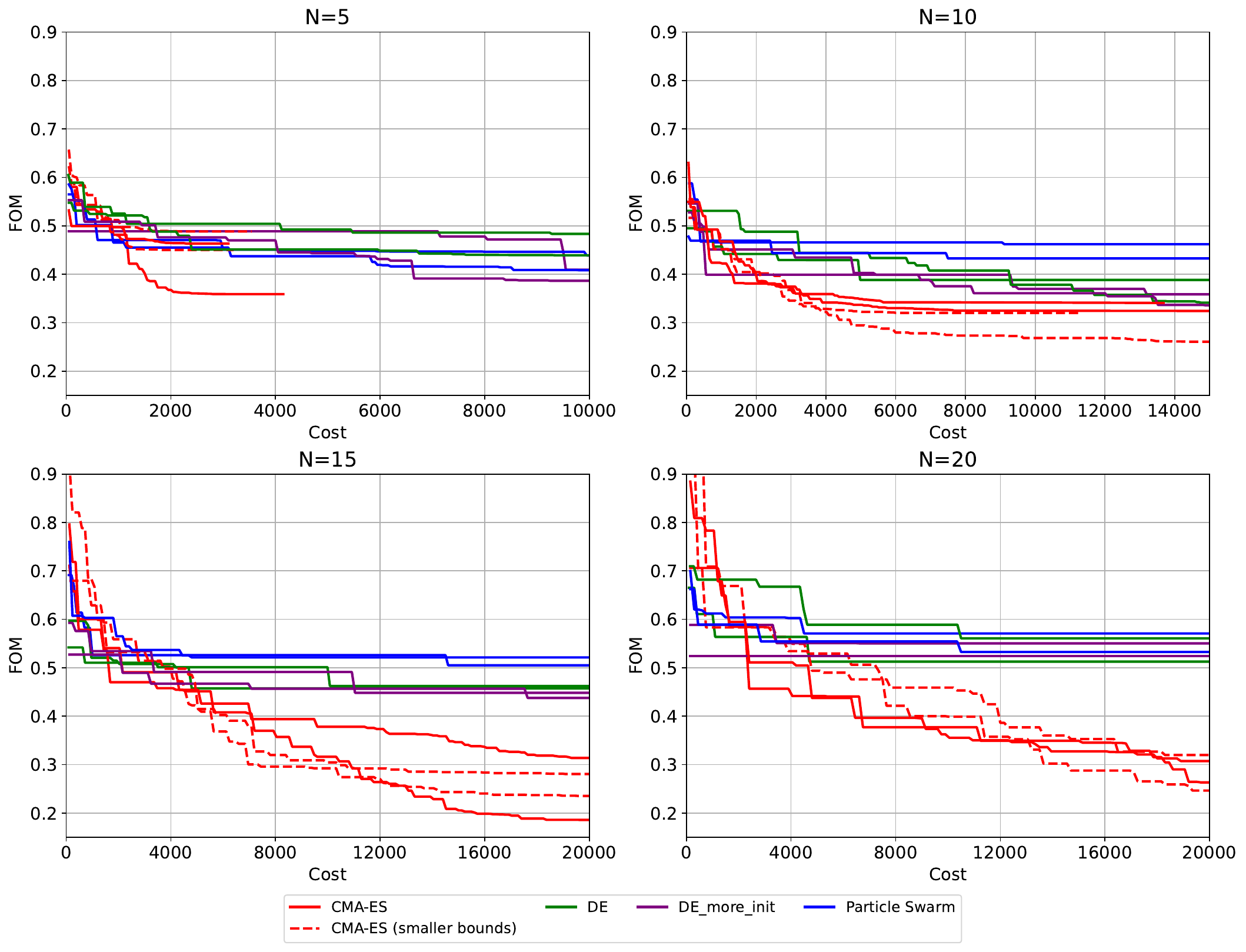}
        \caption{Comparison of the population based methods for different number of cycles ($N$). $N = 5$ (Top-Left), $N = 10$ (Top-Right), $N = 15$ (Bottom-Left), $N = 20$ (Bottom-Right). Note that, for $N=5$ and $N=10$ the CMA-ES algorithm sometimes converges before reaching the computational budget.}\label{fig:comparison_N}
    \end{figure}

    In figure~\ref{fig:comparison_N} we present optimization results for $N=5$, $N=10$, $N=15$, and $N=20$ for the population based optimization algorithms (the different variants of DE, CMA-ES, and PSO). Each algorithm is run twice for different initial populations. For all configurations CMA-ES generally converges the fastest, similar to what has been observed in section~\ref{sec:smooth-beam-results}. We note, however, that for $N=15$ and $N=20$ (i.e.~a higher-dimensional optimization problem) CMA-ES is the only algorithm that, within the computational budget considered, performs satisfactorily.

    \begin{figure}[H]
        \centering
        \includegraphics[width = 0.9\textwidth]{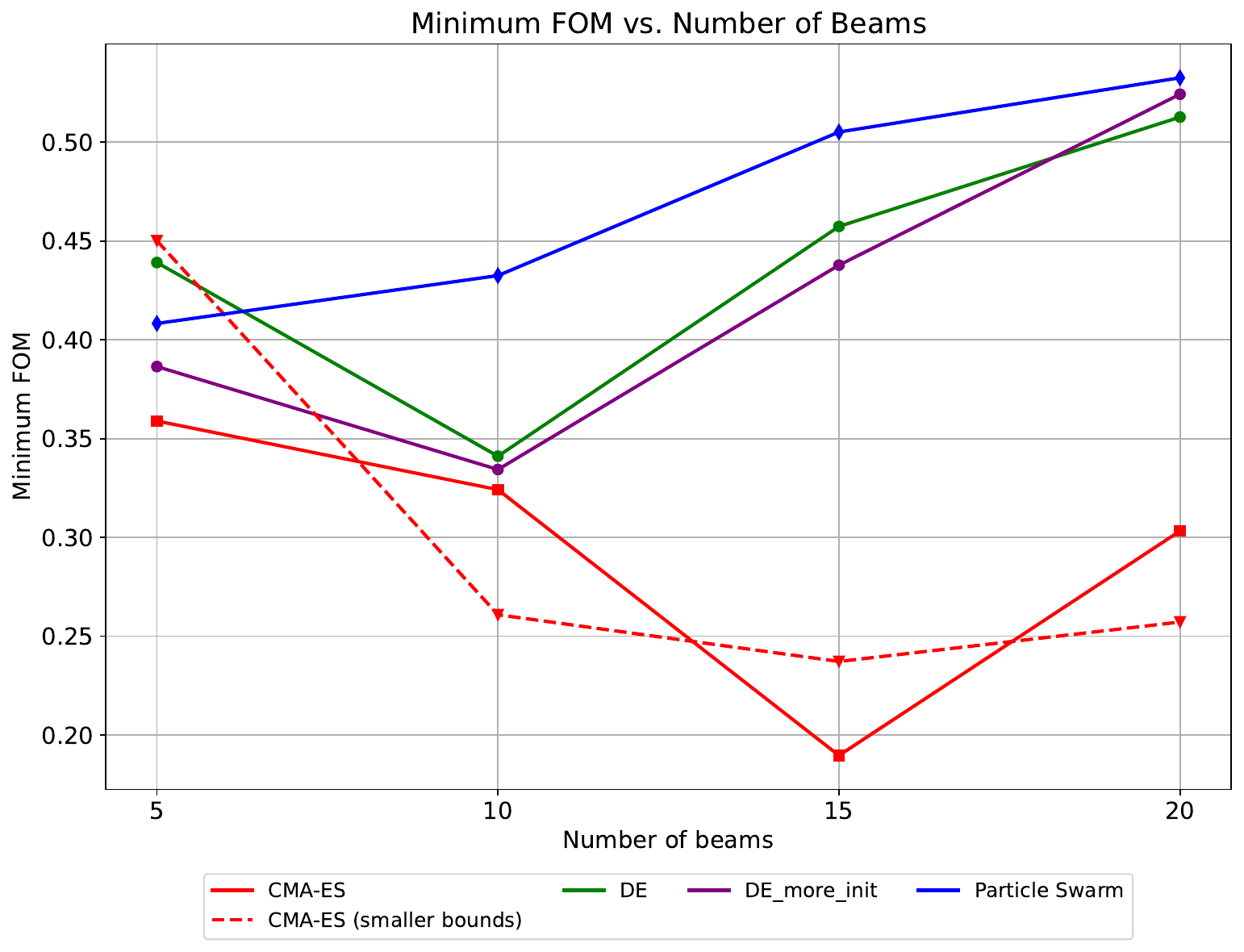}
        \caption{Comparison of the best FOM values found for different numbers of cycles (N). For $N=5$ we compare the performance at $10000$ cost, for $N = 10$ at $15000$ and for $N=15, N=20$ at $20000$ cost (or at the point of convergence).}\label{fig:min_fom_vs_beams}
    \end{figure}

In figure~\ref{fig:min_fom_vs_beams} we compare the best FOM found by the analyzed algorithms for different $N$. For $N = 5$ the computational budget is set to $10000$ (goal) function evaluations, for $N = 10$ it is set to $15000$ function evaluations and for $N=15, N=20$ it is set to $20000$ function evaluations. We choose larger budgets for the higher dimensional problems, since they are more difficult to optimize. In figure~\ref{fig:comparison_N} we can see that for $N=5$ and $N=10$, CMA-ES converges before reaching the computational budget. In that case we compare the FOM of CMA-ES at the point of convergence. The FOM obtained for CMA-ES is monotonically decreasing (i.e.~improves) in $N$ until $N = 15$ where it reaches a minimum of $0.189$. After that the performance decreases again for $N = 20$. We conjecture that this is due to the fact that the optimization problem becomes more difficult with increasing $N$. That is, no good parameter set can be found within the computational budget for $N=20$. An alternative explanation is that the optimal number of cycles is around $N = 15$. However, since we see a similar behavior for DE earlier (going from $N = 10$ to $N=15$) and PSO does not perform satisfactorily even for $N=10$ we tend to lean towards the former hypothesis.

CMA-ES clearly outperforms DE and PSO here. Both DE variants perform generally better than PSO. We reach the best performance of DE\_more\_init at $N = 10$ with a FOM of $0.33$. The DE variants perform similarly to what we observed in section~\ref{sec:smooth-beam-results}. DE gets outperformed by DE\_more\_init for $N<20$. CMA-ES and CMA-ES with smaller bounds both perform well. It is not surprising that there is no clear winner between the two approaches. On the one hand, the smaller bounds are more restrictive and therefore limit the search space. On the other hand, the reduced search space makes the optimization problem easier to solve. It is also worth noting that CMA-ES is a stochastic algorithm, and therefore we expect the results to vary considerably between different runs.
    \begin{figure}[H]
        \centering
        \includegraphics[width = 1.0\textwidth]{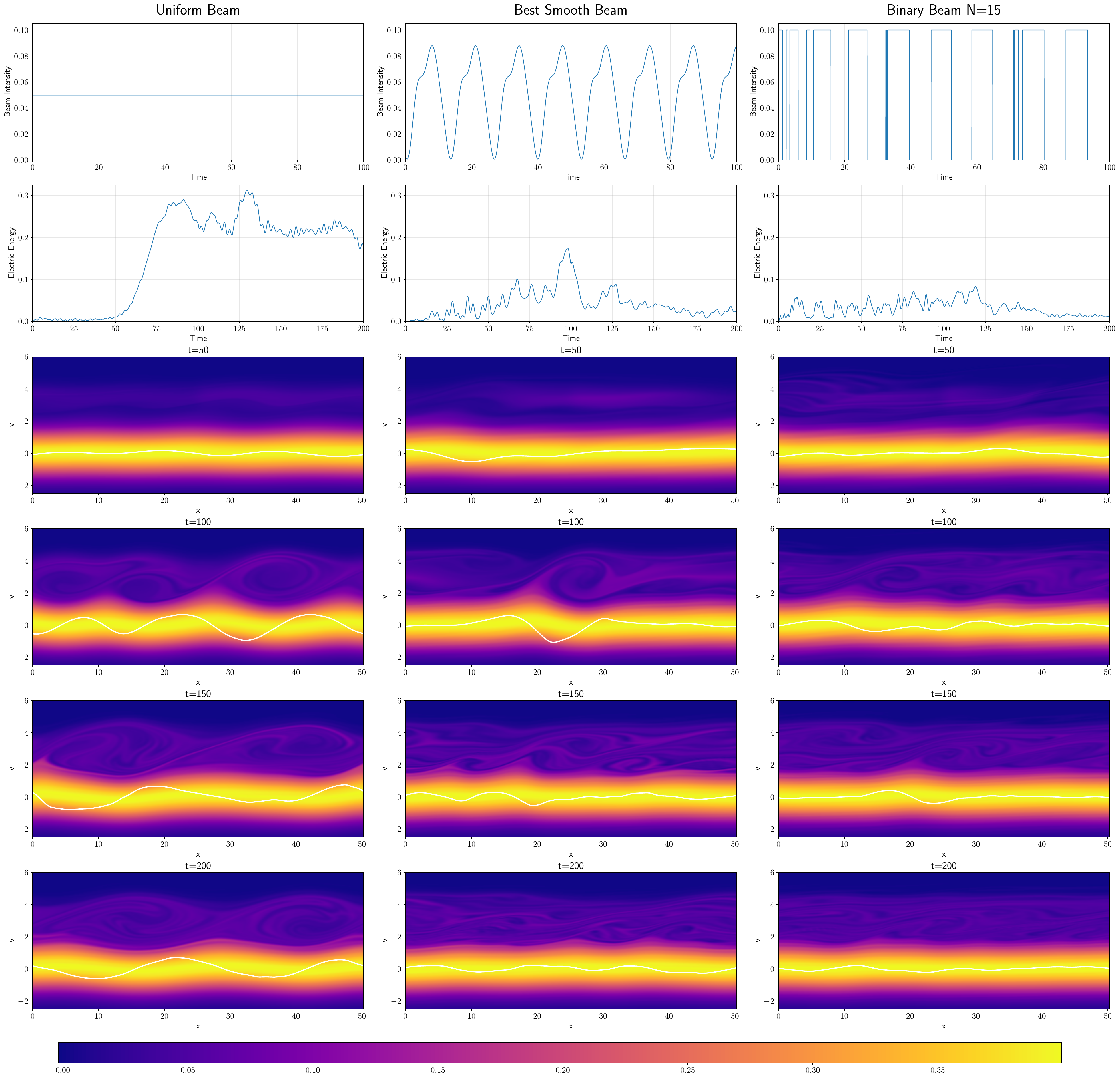}
        \caption{Comparison of beam profiles, electric energies and plasma distributions for times $t = 50, 100, 150, 200$. From left to right we compare the uniform beam ($FOM = 1$), the best smooth beam (found by DE, $FOM = 0.285$), and the best binary beam with $N=15$ (found by CMA-ES, $FOM = 0.189$). All the parameters are listed in table~\ref{tab:2}.}\label{fig:bi_ee_plasma2}
    \end{figure}
In figure~\ref{fig:bi_ee_plasma2} we compare the best binary beam to the uniform and the best smooth beam. This clearly shows the improvement attained by the binary beam, which has a FOM of $0.189$ compared to the best smooth beam which has a FOM of $0.285$ (see section~\ref{sec:smooth-beam-results} and~\cite{10.1063/5.0238408}). This shows that it is possible to obtain excellent performance of our beam modulation by simply switching the beam ON and OFF.
    \begin{figure}[H]
        \centering
        \includegraphics[width = 1.0\textwidth]{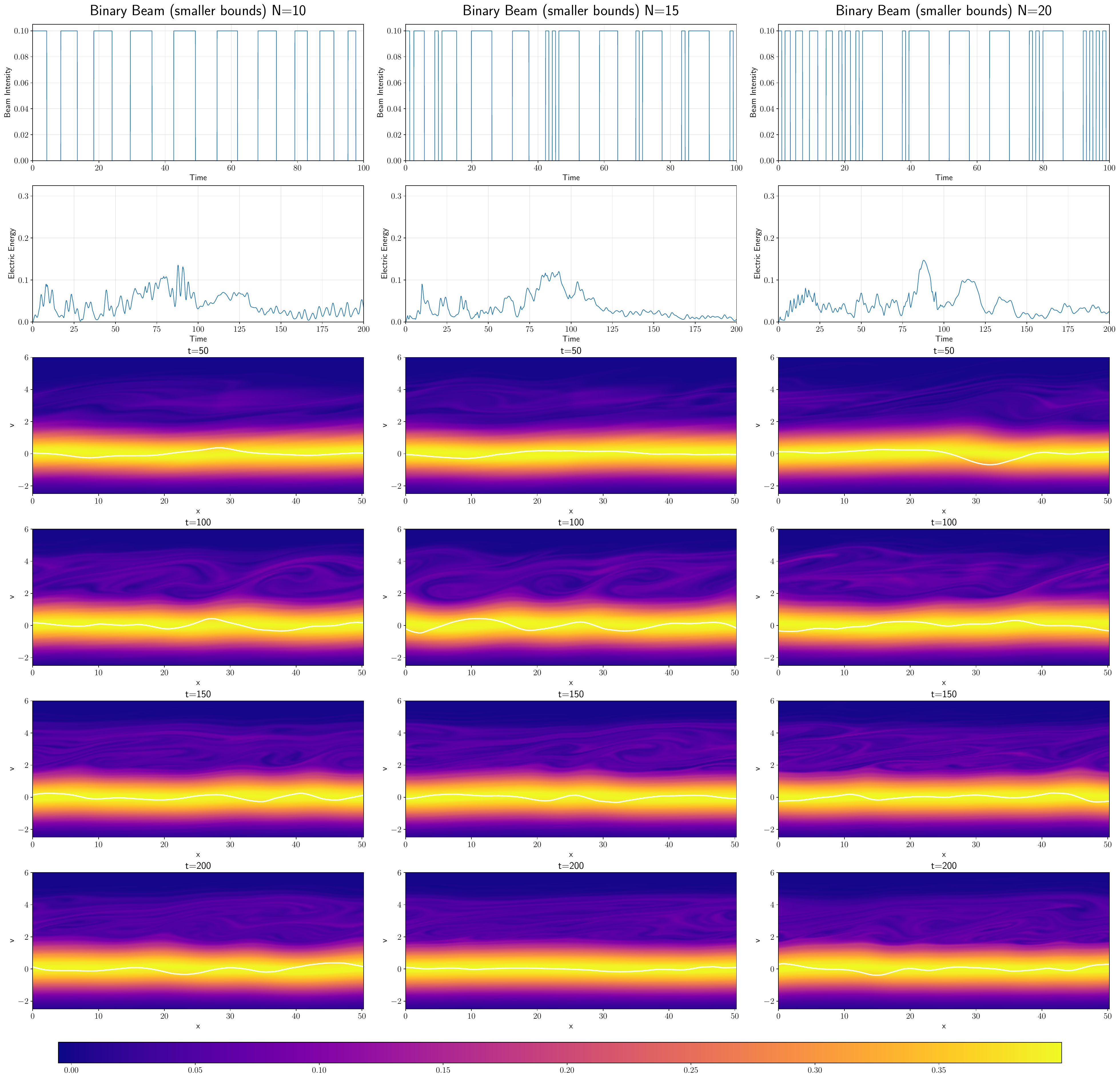}
        \caption{Comparison of the best binary beam profiles with smaller bounds, electric energies and plasma distributions for times $t = 50, 100, 150, 200$ for different numbers of cycles $N$. All of these configurations were found by CMA-ES. From left to right we compare the binary beam with $N=10$ ($FOM = 0.261$), $N=15$ ($FOM = 0.235$) and the binary beam with $N=20$ ($FOM = 0.261$). All the parameters are listed in table~\ref{tab:2}.}\label{fig:bi_ee_plasma3}
    \end{figure}
The only disadvantage of the found binary beams is that some cycles have a rather short duration. The shortest cycle of the binary beam profile in figure~\ref{fig:bi_ee_plasma2} has a duration of $0.26$. This can be a challenge to realize experimentally. Thus, we show the best binary beam profiles with smaller bounds (i.e.~$\epsilon=\tfrac{4}{25}$) for different numbers of cycles in figure~\ref{fig:bi_ee_plasma3}. All of these configurations were found by CMA-ES.~The best performance with a FOM of $0.235$ is achieved with $15$ cycles. The duration of the cycles, as can be seen from figure~\ref{fig:bi_ee_plasma3}, is significantly longer than the cycles of the binary beam in figure~\ref{fig:bi_ee_plasma2}. The duration of the shortest cycle in this case is $1.99$, so more than $7$-times longer than the shortest cycle of the beam in~\ref{fig:bi_ee_plasma2}. We also highlight the good performance of the $10$-cycle case, which has a better FOM ($0.261$) than the best smooth beam ($0.285$) and very long cycle durations with the shortest cycle having a duration of $4.74$. Note how the $20$-cycle beam reaches the same performance but with twice as many cycles that are much shorter, with the shortest cycle having a duration of $1.95$.

    \begin{figure}[H]
        \centering
        \includegraphics[width = 1.0\textwidth]{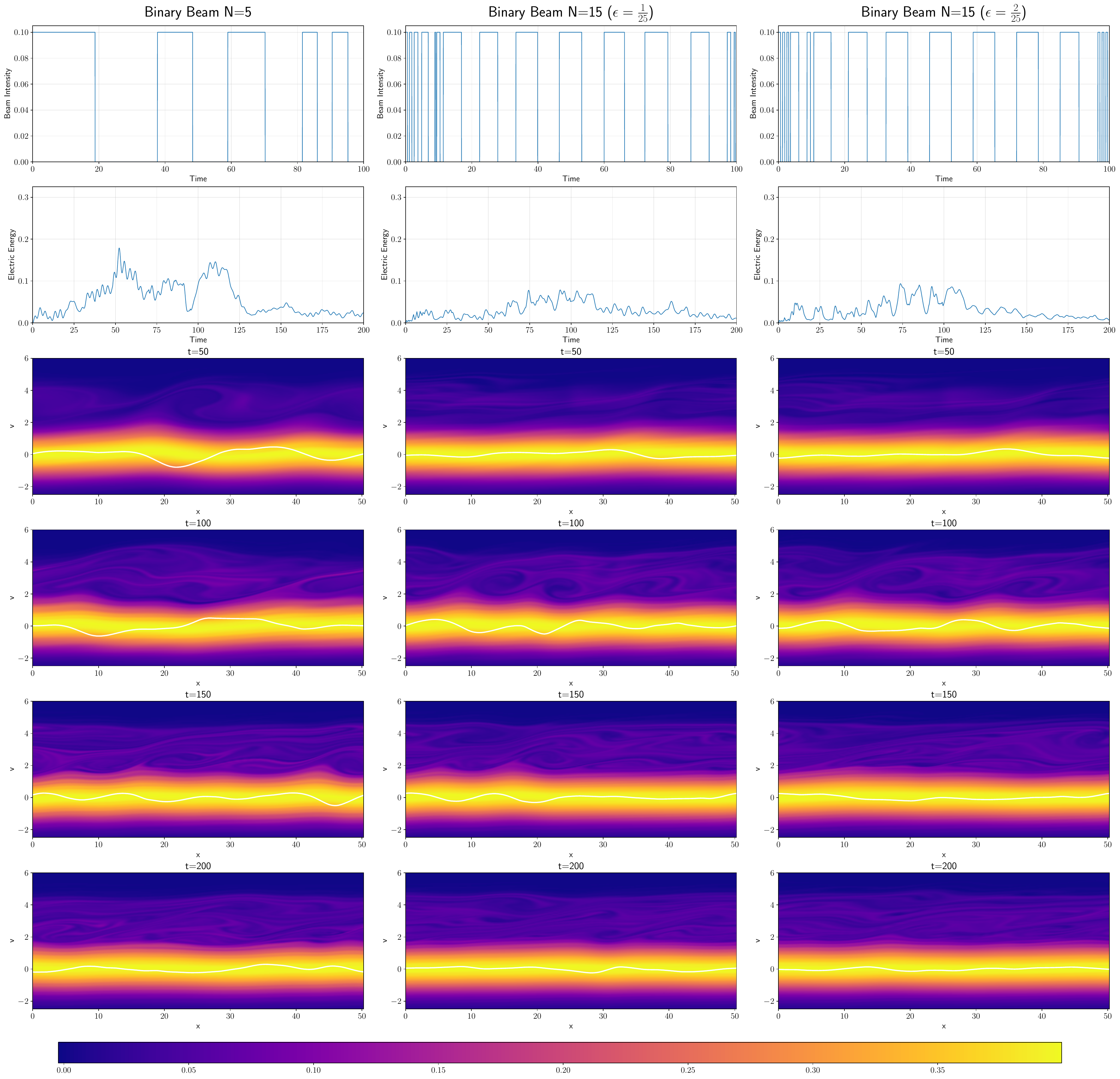}
        \caption{Comparison of more binary beam profiles, electric energies and plasma distributions for times $t = 50, 100, 150, 200$. From left to right we compare the binary beam with $N=5$ ($FOM = 0.359$) with two binary beams with $N=15$. For the two beams in the second and third columns we used smaller $\epsilon$ as lower bounds, hence allowing shorter cycles. We used $\epsilon = \frac{1}{25}$ and $\epsilon = \frac{2}{25}$ leading to $FOM = 0.183$ and $FOM = 0.197$, respectively. All the parameters are listed in table~\ref{tab:2}.}\label{fig:bi_ee_plasma4}
    \end{figure}
    
In figure~\ref{fig:bi_ee_plasma4} we compare additional binary beams found by the optimizer during our numerical experiments. From left to right we compare the binary beam with $N=5$, reaching $FOM = 0.359$, with two binary beams with $N=15$. However, for these two beams we used smaller $\epsilon$ as the smaller bounds case, thus, shorter cycles were allowed. For the beam in the middle we used $\epsilon = \frac{1}{25}$ and this beam achieved the best overall FOM of $0.183$, whereas the shortest cycle in this case is of duration $0.60$. For the beam on the right we used $\epsilon = \frac{2}{25}$ and this configuration resulted in a FOM of $0.197$. This beam has a significantly longer shortest cycle of $1.08$. The two $15$-cycle beams are mentioned here due to their excellent performance despite having longer cycle durations than the binary beam in~\ref{fig:bi_ee_plasma2}. The cycle durations of the $5$-cycle case are obviously very long, with the shortest cycle being of length $8.99$, making it attractive for experiments. However, this beam has a significantly worse FOM compared to the others.~Nevertheless, it is still interesting to note that it reduces the severity of the instability by almost two-thirds. This tells us that even in a relatively low-dimensional parameter space interesting beam profiles can be found.

\section{Conclusion}
    We have investigated the optimization of beam profiles for a Vlasov-Poisson plasma model and compared several derivative-free and gradient-based methods for two different parametrizations. The resulting objective landscapes are highly non-convex and contain many shallow and deep local minima. This strongly limits the reliability of local methods such as L-BFGS-B and even Basinhopping, which typically converge to nearby local minima rather than to the best attainable solution. Among the optimization methods tested, CMA-ES and Differential Evolution yield the best performance. CMA-ES, in particular, performs best for problems with either a more limited computational budget or for higher-dimensional problems. Differential Evolution, in some situations, can find deeper minima than CMA-ES.

    We have considered binary beams, which could be easier to realize experimentally, but lead to a higher-dimensional optimization problem. We have found such beam profiles (using CMA-ES and Differential Evolution) that reduce the severity of the bump-on-tail instability by more than a factor of five, which significantly improves upon what has been reported in the literature.

    We highlight that this is a difficult optimization problem and the methods employed require at least multiple thousands of evaluations of the goal function (i.e.~simulations of the Vlasov equation). Thus, in order to tackle higher-dimensional plasma dynamics (i.e.~more than a single dimension in space and velocity), we will investigate complexity and model order reduction techniques in the future. Also, it would be of interest, in future work, to study models that are closer to the neutral beam injection performed in tokamaks. This, in particular, requires the multi-scale modeling of two species (electrons and ions) as well as considering the dynamics of the plasma in strong magnetic fields.

\section{Appendix}
    \begin{table}[H]
        \centering
        \caption{Parameters for the beam profiles shown in Figures~\ref{fig:bi_ee_plasma1},~\ref{fig:bi_ee_plasma3},~\ref{fig:bi_ee_plasma2} and~\ref{fig:bi_ee_plasma4}.}\label{tab:2}
        \small
        \setlength{\tabcolsep}{3pt}
        \renewcommand{\arraystretch}{1.0}
        \newcommand{\tval}[2]{\makebox[0.27\linewidth][l]{$t_{#1} = #2$}}
        \newcommand{\pval}[1]{\makebox[0.27\linewidth][l]{#1}}
        \begin{tabularx}{\textwidth}{@{}p{3.5cm}X@{}}
            \hline
            \textbf{Sine} & \pval{$\omega = 0.48221$}\pval{$b_1 = 0.03$}\\
            \hline
            \textbf{Optimized} & \pval{$\omega = 0.47789$}\pval{$a_1 = -0.037118$}\pval{$a_2 = -0.010956$}\\
            & \pval{$b_1 = -0.015525$}\pval{$b_2 = -0.00095985$}\pval{$b_3 = -0.0065273$}\\
            \hline
            \textbf{Best Smooth} & \pval{$\omega = 0.47734131$}\pval{$a_1 = -0.03684612$}\pval{$a_2 = -0.01002340$}\\
            & \pval{$b_1 = -0.01000293$}\pval{$b_2 = -0.00123800$}\pval{$b_3 = -0.00430074$}\\
            \hline
            \textbf{Binary Beam} $N = 5$ & \tval{1}{0.99495705}\tval{2}{0.56127325}\tval{3}{0.59486417}\\
            & \tval{4}{0.23725157}\tval{5}{0.25122083}\\
            \hline
            \textbf{Binary Beam} $N=10$ & \tval{1}{0.64796088}\tval{2}{0.76035923}\tval{3}{0.84045809}\\
            \textbf{(smaller bounds)} & \tval{4}{0.99998475}\tval{5}{0.99907897}\tval{6}{0.94173459}\\
            & \tval{7}{0.84988764}\tval{8}{0.57655125}\tval{9}{0.64674143}\\
            & \tval{10}{0.36129834}\\
            \hline
            \textbf{Binary Beam} $N=15$ & \tval{1}{0.17620656}\tval{2}{0.08311684}\tval{3}{0.38616139}\\
            & \tval{4}{0.15727158}\tval{5}{0.80130459}\tval{6}{0.85722722}\\
            & \tval{7}{0.02126564}\tval{8}{0.02002059}\tval{9}{0.99935526}\\
            & \tval{10}{0.93356903}\tval{11}{0.95622041}\tval{12}{0.02075461}\\
            & \tval{13}{0.18159037}\tval{14}{0.99853134}\tval{15}{0.99978676}\\
            \hline
            \textbf{Binary Beam} $N=15$ & \tval{1}{0.20047582}\tval{2}{0.50777213}\tval{3}{0.17393597}\\
            \textbf{(smaller bounds)} & \tval{4}{0.71261847}\tval{5}{0.99455550}\tval{6}{0.80455749}\\
            & \tval{7}{0.16225354}\tval{8}{0.16000333}\tval{9}{0.98577909}\\
            & \tval{10}{0.88230673}\tval{11}{0.16185893}\tval{12}{0.94586505}\\
            & \tval{13}{0.17228503}\tval{14}{0.99948132}\tval{15}{0.16204928}\\
            \hline
            \textbf{Binary Beam} $N=15$ & \tval{1}{0.08511656}\tval{2}{0.10191165}\tval{3}{0.16300686}\\
            \textbf{($\epsilon = \frac{1}{25}$)} & \tval{4}{0.28476738}\tval{5}{0.04283528}\tval{6}{0.14223300}\\
            & \tval{7}{0.79023544}\tval{8}{0.78718796}\tval{9}{0.94380783}\\
            & \tval{10}{0.97532609}\tval{11}{0.88563896}\tval{12}{0.99976859}\\
            & \tval{13}{0.79043962}\tval{14}{0.14772989}\tval{15}{0.05408985}\\
            \hline
            \textbf{Binary Beam} $N=15$ & \tval{1}{0.09404389}\tval{2}{0.09817993}\tval{3}{0.08279646}\\
            \textbf{($\epsilon = \frac{2}{25}$)} & \tval{4}{0.37918971}\tval{5}{0.15433245}\tval{6}{0.79219356}\\
            & \tval{7}{0.86266179}\tval{8}{0.99852147}\tval{9}{0.99751010}\\
            & \tval{10}{0.99767369}\tval{11}{0.99556572}\tval{12}{0.86497764}\\
            & \tval{13}{0.09505275}\tval{14}{0.08935779}\tval{15}{0.08193800}\\
            \hline
            \textbf{Binary Beam $N=20$} & \tval{1}{0.16163761}\tval{2}{0.26412967}\tval{3}{0.34149637}\\
            \textbf{(smaller bounds)} & \tval{4}{0.41304320}\tval{5}{0.31449544}\tval{6}{0.16030328}\\
            & \tval{7}{0.26150926}\tval{8}{0.16836974}\tval{9}{0.98630377}\\
            & \tval{10}{0.16391255}\tval{11}{0.99977690}\tval{12}{0.99984740}\\
            & \tval{13}{0.98519022}\tval{14}{0.16066296}\tval{15}{0.17828661}\\
            & \tval{16}{0.99883359}\tval{17}{0.16053184}\tval{18}{0.16107889}\\
            & \tval{19}{0.16008810}\tval{20}{0.17014673}\\
            \hline
        \end{tabularx}
    \end{table}
\printbibliography{}
\end{document}